\documentclass[a4paper,11pt,reqno]{amsart}
\usepackage[utf8]{inputenc}
\usepackage{amsmath,amssymb,amsthm} 
\usepackage{graphics} 

\usepackage[colorlinks,citecolor=blue]{hyperref}

\usepackage{epsfig}
\usepackage{color}
\usepackage{esint}

\numberwithin{equation}{section}

\newtheorem{theorem}{Theorem}[section]
\newtheorem{lemma}[theorem]{Lemma}
\newtheorem{proposition}[theorem]{Proposition}
\newtheorem{corollary}[theorem]{Corollary}

\theoremstyle{definition}
\newtheorem{definition}[theorem]{Definition}

\newcommand{\R}{\mathbb{R}}
\newcommand{\N}{\mathbb{N}}

\newcommand{\Z}{\mathbb{Z}}

\newcommand{\eps}{\varepsilon}

\newcommand{\rst}{\llcorner}

\newcommand{\be}{\begin{equation}}
\newcommand{\ee}{\end{equation}}

\DeclareMathOperator{\supp}{supp}

\newcommand\lt{\left}
\newcommand\rt{\right}

\def\les{\lesssim}
\def\ges{\gtrsim}

\def\Id{{\rm Id}}

\def\EE{\mathbb{E}}
\def\T{\mathbb{T}}

\def\one{\mathbf{1}}
\def\PP{\mathbb{P}}

\def\hE{\widehat{E}}

\newcommand{\step}[1]{\medskip\noindent\textbf{Step #1. }}

\def\spt{ {\rm Spt}\,}
\def\Var{ {\rm Var}}

\newcommand{\Dat}{\mathfrak{D}}          

\title[Large-scale regularity for optimal $p$-transport maps]{Large-scale regularity theory for optimal $p-$transport maps and application to concentration in the matching problem}

\author{Michael Goldman}
\address{Univ Paris Est Creteil, CNRS, LAMA, F-94010 Creteil, France}
\email{michael.goldman@cnrs.fr}

\author{Lukas Koch}
\address{Department of Mathematics, University of Sussex, Brighton, BN1 9QH, United Kingdom}
\email{lukas.koch@sussex.ac.uk}

\author{Dario Trevisan}
\address{Dipartimento di Matematica, Universit\`a di Pisa, Largo Bruno Pontecorvo 5, 56127 Pisa, Italy}
\email{dario.trevisan@unipi.it}

\begin{document}

\begin{abstract}
We consider the random bipartite matching problem between two independent families of uniformly distributed points 
in dimension three and higher, for the cost given by the $p$-th power of the distance. 
We show that, up to logarithmic factors, every moment of the displacement of the optimal matching is of the order of the 
typical interparticle distance. This is a higher integrability, or Calder\'on--Zygmund type, estimate for optimal transport, 
and it settles, on the torus and up to logarithmic factors, a recent conjecture of D'Achille, Mattesini and the third author. 
As a consequence, we obtain concentration of the matching cost around its mean for every $p\ge2$. 
The proof combines a deterministic large-scale regularity theory for optimal transport with $p$-costs, 
based on a harmonic approximation by a non-linear elliptic equation of $p$-Laplacian type, with tail bounds on local matching costs.
\end{abstract}
\maketitle

\section{Introduction}\label{sec:intro}

Let $d\ge3$, $L\ge2$ with $n=L^d\in\N$, and let $\lambda=\sum_{i=1}^n\delta_{X_i}$, $\mu=\sum_{j=1}^n\delta_{Y_j}$, where $(X_i)$ and $(Y_j)$ are independent families of i.i.d.\ uniform points on the flat torus $\T_L=(\R/L\Z)^d$. For $p\ge1$ we consider the random bipartite matching problem
\[
F_{L,p}=\frac1{L^d}W_{p,\T_L}^p(\lambda,\mu)=\min_{\sigma}\frac1n\sum_{i=1}^nd_{\T_L}(X_i,Y_{\sigma(i)})^p,
\]
where the minimum runs over permutations of $\{1,\dots,n\}$, and we denote by $\pi_p$ the corresponding optimal coupling.
 While it has been known for a long time that $\EE\lt[F_{L,p}\rt]\simeq 1$ 
and even that $\EE\lt[F_{L,p}\rt]$ converges to a finite constant as $L\to\infty$, see \cite{barthe2013combinatorial,dereich2013constructive,goldman2021convergence}, 
the fluctuations of $F_{L,p}$  around its mean are less understood, see however \cite{feng2026fluctuations} for recent results in the case $p=2$. 
Our main result is a concentration estimate for $F_{L,p}$ which holds for every $p\ge2$.

\begin{theorem}\label{thm:concentration}
Let $d\ge3$ and $p\ge2$. There is $C=C(d,p)$ such that, for $L\ge2$ with $L^d\in\N$,
\[
\Var(F_{L,p})\le C L^{2-d}(\log L)^{2(p-1)(1+\frac1d)}.
\]
In particular, $\PP\big[|F_{L,p}-\EE F_{L,p}|\ge t\big]\le Ct^{-2}L^{2-d}(\log L)^{2(p-1)(1+\frac1d)}$ for every $t>0$.
\end{theorem}

Concentration for the matching cost has been known for a while for $p<d$, as well as for $p=d=2$, see \cite{AmStTr16} and also \cite{goldman2021convergence}. The difficulty for larger $p$ is that the natural route to concentration, through the Poincar\'e (or Efron--Stein) inequality, requires a bound on the $(2p-2)$-th moment of the displacement of the $p$-optimal matching, which is not provided by the optimality of $\pi_p$. For general $p$, only tail bounds for the matching cost were recently obtained by Huesmann, Mattesini and Otto \cite{huesmann2024there} (let us stress that they were mostly interested in the case $d=2$ which is much more delicate).
 We give simpler proofs of these tail bounds, see the discussion below. We restrict ourselves  in this paper to $p\ge2$ since this is the only range in which Theorem \ref{thm:concentration} is new.

Our approach builds on the ideas of D'Achille, Mattesini and the third author \cite{d2026concentration}, who reduced concentration to a $p\to q$ 
transfer principle. Our second main result proves their conjecture for $p\ge2$, on the torus and up to logarithmic factors.

\begin{theorem}\label{thm:main}
Let $d\ge3$, $p\ge2$ and $q\ge1$. There is $C=C(d,p,q)$ such that, for $L\ge2$ with $L^d\in\N$,
\begin{equation}\label{mainclaim}
\EE\lt[ \frac{1}{L^d} \int_{\T_L\times \T_L} d_{\T_L}(x,y)^q\, d \pi_p\rt]^{1/q}\le C(\log L)^{1+\frac1d}.
\end{equation}
\end{theorem}
 Theorem \ref{thm:concentration} follows from Theorem \ref{thm:main} with $q=2p-2$ as in \cite{d2026concentration}. We provide a short proof for the reader's convenience. 
  Theorem \ref{thm:main} can be seen as a higher integrability, or Calder\'on--Zygmund type, estimate for optimal transport.

To prove Theorem \ref{thm:main} we extend the strategy of the first author and Huesmann \cite{goldman2022fluctuation}, developed there for $p=2$. By stationarity, it is enough to bound the displacement of the points in a ball of logarithmic size around the origin. This is obtained by combining a deterministic large-scale regularity theory for $p$-optimal transport,
 in the spirit of the first half of \cite{goldman2021quantitative}, with a stochastic input. To state the deterministic result, let $Q_R=(-R,R)^d$, $\#_R=(B_R\times\R^d)\cup(\R^d\times B_R)$, and for a measure $\nu$ let $\kappa_{\nu,Q_R}=\nu(Q_R)/|Q_R|$ and
\[
\Dat^\square_R(\nu)=\frac1{R^d}W_p^p\big(\nu\rst Q_R,\kappa_{\nu,Q_R}\,dx\rst Q_R\big)+R^p\,|\kappa_{\nu,Q_R}-1|^p .
\]

\begin{theorem}\label{thm:iteration}
Let $d\ge3$, $p\ge2$ and $C_D,C_E\ge1$. There exist $1\le M_0,C_\star,C_{\rm exc}<\infty$, depending only on $d,p,C_D,C_E$, such that the following holds.
 Let $r_\ast\ge2$, $R_0\ge M_0r_\ast$, let $\lambda,\mu$ be measures on $\R^d$, and let $\pi$ be an optimal coupling between $\lambda$ and $\mu$
  for the cost $|y-x|^p$. Assume:
\begin{enumerate}
\item[(H1)] for every dyadic $R=2^j\in[r_\ast,20R_0]$, $j\in\Z$, $\Dat^\square_R(\lambda)+\Dat^\square_R(\mu)\le C_Dr_\ast^p$;
\item[(H2)] $\displaystyle\frac1{(4R_0)^d}\int_{\#_{4R_0}}|y-x|^p\,d\pi\le C_Er_\ast^p$.
\end{enumerate}
Then, with $R_\star=C_\star^{-1}r_\ast\log(R_0/r_\ast)$,
\begin{equation}\label{eq:excessdecay}
\sup_{R\in[R_\star,R_0]}\ \inf_{b\in \R^d}\ \frac{1}{\log(4R_0/R)^{1-\frac2p}}\Big(\frac1{R^d}\int_{\#_R}|y-x-b|^p\,d\pi\Big)^{\frac1p}\le C_{\rm exc}\,r_\ast,
\end{equation}
and
\begin{equation}\label{eq:finaldisplacement}
\sup\big\{|y-x|\colon (x,y)\in\spt\pi\cap\#_{R_\star}\big\}\le C_\star\,r_\ast\log(R_0/r_\ast).
\end{equation}
\end{theorem}

Theorem \ref{thm:iteration} is a large-scale regularity result. Hypothesis (H1) says that $\lambda$ and $\mu$ are close to the Lebesgue measure at all dyadic scales between $r_\ast$ and $R_0$, and (H2) that the transport energy is small at the largest scale. The estimate \eqref{eq:excessdecay} then states that, at every scale above $R_\star$, the displacement $y-x$ is close in an averaged sense to a constant vector $b_R$, and \eqref{eq:finaldisplacement} follows from an $L^\infty$ bound at scale $R_\star$. The constant $b_R$ may drift by $O(r_\ast)$ at each dyadic scale, which is the origin of the logarithmic loss in \eqref{eq:finaldisplacement}. This large-scale regularity theory is inspired by quantitative stochastic homogenization, see \cite{armstrong2019quantitative}. In this terminology, Theorem \ref{thm:iteration} is a large-scale $C^{0,\alpha}$- theory for the displacement. To obtain a $C^{1,\alpha}$-theory, as in the small-scale regularity theory of \cite{GO}, one would also need to take into account the affine part of the displacement, that is its first-order Taylor expansion. For non-quadratic costs this leads to severe difficulties, since $p$-costs do not share the invariance of the quadratic cost under affine changes of coordinates.

Theorem \ref{thm:iteration} is obtained by iterating a one-step improvement (Proposition \ref{onestep}), which relies on a harmonic approximation result (Theorem \ref{harmapprox}). After a change of coordinates factoring out the constant displacement $b$, the problem at each scale becomes an optimal transport problem for the shifted cost $c_b$ of \eqref{eq:defhatc}, which behaves quadratically for $|z|\ll|b|$ and like $|z|^p$ for $|z|\gg|b|$. Theorem \ref{harmapprox} extends the harmonic approximation of the second author \cite{koch2024geometric} for $p$-costs to this setting. While the proof is parallel to \cite{koch2024geometric} (see also \cite{GO,goldman2021quantitative,KoOt} for the quadratic case), making all the estimates uniform in $b$ requires some care, and in particular a heavy use of the $V$-functions from the regularity theory of $p$-Laplacian type equations \cite{GiaquintaModica1986,AcerbiFusco1989}. 
The harmonic approximation is given by a shifted $p'$-Laplace equation, and we need regularity estimates for it which are uniform in $b$. As these were not available in exactly this form, we derive them from the literature in Appendix \ref{app:PDE}. 
 We expect the  harmonic approximation result (Theorem \ref{harmapprox}) to also hold for $p\in(1,2)$. However, as for the $p$-Laplacian, the sub-quadratic case requires
 slightly different computations, see for instance \cite{duzaar2004} where both cases are treated, and we did not include it.

The stochastic input consists of tail bounds for the local matching costs $\Dat^\square_R(\lambda)$ and $\Dat^\square_R(\mu)$, which ensure that the assumptions of Theorem \ref{thm:iteration} hold with $r_\ast\sim(\log L)^{1/d}$ with high probability. Tail bounds of this type were proven in \cite{huesmann2024there}, where the main case of interest was $d=2$, which cannot be handled by our techniques. In Appendix \ref{app:moments} we give much simpler proofs for $d\ge3$. As in \cite{goldman2021convergence}, we use a multiscale decomposition into dyadic cubes, which reduces the proof to concentration estimates for binomial random variables. If one stops at logarithmically large scales, which is enough for our application since the iteration stops at such scales anyway, a union bound and Bernstein's inequality suffice (Theorem \ref{thm:momentlog}). At these large scales, our bounds are similar to the ones of Garc\'ia Trillos and Slep\v{c}ev \cite{Garcia-Slepcev} for the $\infty$-transport distance. Going down to unit scale (Theorem \ref{thm:momentunit}) is slightly more delicate. Since a union bound no longer suffices, we exploit more cancellations: at small scales many cubes are involved, and the sum of their contributions exhibits a CLT type behavior. We obtain stretched exponential tail bounds, improving on the polynomial ones of \cite{huesmann2024there}.

The paper is organized as follows. In Section \ref{sec:prelim} we introduce the shifted costs and collect preliminaries on optimal transport. Section \ref{sec:deterministic} contains the harmonic approximation, the one-step improvement and the proof of Theorem \ref{thm:iteration}. Theorems \ref{thm:main} and \ref{thm:concentration} are proven in Sections \ref{sec:higher} and \ref{sec:concentration}. Appendix \ref{app:PDE} contains the PDE estimates, and Appendix \ref{app:moments} the tail bounds.

\section{Preliminaries}\label{sec:prelim}

Throughout, $d\ge3$ and $p\ge2$. We write $p'=p/(p-1)\in(1,2]$. Constants called $C,c$ depend only on $d$ and $p$ unless otherwise indicated, and $A\les B$ means $A\le CB$.
 The symbol $A\simeq B$ means $A\les B\les A$.  We use $A\ll B$ as an hypothesis saying that there exists $\eps>0$ such that if $A\le \eps B$ then the conclusion holds.

\subsection{Notation}
For $R>0$, $B_R(x_0)$ denotes the open ball of radius $R$ around $x_0$. If $x_0=0$, we write $B_R = B_R(0)$. Similarly, $Q_R(x_0)=x_0+(-R,R)^d$ and $Q_R=Q_R(0)$. Denote for $R>0$, $\#_R = (B_R\times \R^d)\cup (\R^d\times B_R)$. If $\lambda$ is a measure admitting a density, we will frequently misuse notation and denote its density by $\lambda$. For a measure $\nu$ and a bounded set $U$ with $|U|>0$ we write
\begin{equation}\label{eq:kappa}
\kappa_{\nu,U}=\frac{\nu(U)}{|U|},\qquad \kappa_{\nu,R}=\kappa_{\nu,B_R}.
\end{equation}

If $\pi$ is a transport plan, we let for $(x,y)\in \supp\,\pi$ and $t\in [0,1]$, $X(t) = (1-t)x+ty$, so that $\dot X=y-x$. We set
\begin{align}\label{def:OmegaR}
\Omega_R = \{(x,y)\colon \exists t\in [0,1], X(t)\in \overline{B_R}\}.
\end{align}
If the radius is clear from the context, we write $\Omega = \Omega_R$. For $(x,y)\in\Omega$, we set
\begin{align*}
t_- = \min_{X(t)\in \overline{B_R}} t,\quad \text{ and } \quad t_+ = \max_{X(t)\in \overline{B_R}} t.
\end{align*}
Since $B_R$ is convex, $X(t)\in\overline{B_R}$ if and only if $t\in[t_-,t_+]$.

\subsection{The cost}\label{sec:prelim-cost}
In this section we introduce the $V$-functions, the shifted cost $c_b$ and its conjugate $c_b^\ast$.
We list the elementary properties of these functions that will be used. They are standard in the regularity theory of $p$-Laplacian type problems, 
where $V_b^2$ and $V_{b,\ast}^2$ appear,  see e.g. \cite{GiaquintaModica1986,AcerbiFusco1989,duzaar2004,DieningEttwein2008,DieningKaplickySchwarzacher2012}; for their use for optimal transport with $p$-costs see \cite{koch2024geometric}.
\subsubsection{The $V$-functions}
It will be convenient to denote for $b\in\R^d$,
\begin{align*}
V_b^2(z) &= (|b|^2+|z|^2)^\frac{p-2} 2|z|^2 \simeq (|b|+|z|)^{p-2}|z|^2 \simeq |z|^p+|b|^{p-2}|z|^2,\\
V_{b,\ast}^2(y) &= (|b|^{p-1}+|y|^2)^\frac{p^\prime-2} 2|y|^2\simeq (|b|^{p-1}+|y|)^{p^\prime-2}|y|^2.
\end{align*}
For $b=0$ we recover $V_0^2(z)=|z|^p$ and $V_{0,\ast}^2(y)=|y|^{p'}$. We record here some elementary properties of these functions without proof. Proofs can be found in the previously cited papers. Both functions are radial, convex and increasing in the modulus. 
 We shall use the following elementary scaling properties: for $\ell\ge 1$,
\begin{equation}\label{eq:scalingV}
V_b^2(\ell z)\le \ell^p V_b^2(z),\qquad V_{b,\ast}^2(\ell z)\le \ell^2 V_{b,\ast}^2(z),\qquad V_{b,\ast}^2(\ell^{-1} z)\le \ell^{-p'} V_{b,\ast}^2(z),
\end{equation}
and the joint homogeneity for $t>0$,
\begin{equation}\label{eq:homogeneityV}
V_{tb}^2(tz)=t^pV_b^2(z),\qquad V_{tb,\ast}^2(t^{p-1}y)=t^pV_{b,\ast}^2(y).
\end{equation}
In particular 
\begin{align}\label{eq:triangleV}
V_b^2(x+y)\lesssim V_b^2(x)+V_b^2(y).
\end{align}
The following Young-type inequalities hold: there is $C>0$ such that for every $\eps\in(0,1]$, $z,w\in\R^d$ and $s,t\ge0$,
\begin{align}\label{eq:YoungsV}
\langle z,w\rangle &\leq \eps V_b^2(z)+C\eps^{-1} V_{b,\ast}^2(w),\qquad
\langle z,w\rangle \leq \eps V_{b,\ast}^2(w)+C\eps^{1-p} V_{b}^2(z),\\
(|b|^{p-1}+t)^{p'-2}ts &\leq (|b|^{p-1}+t)^{p'-2}t^2+(|b|^{p-1}+s)^{p'-2}s^2. \label{eq:YoungsV2}
\end{align}

We also note the elementary inequality, 
\begin{equation}\label{eq:elementary}
V_{b,\ast}^2(|b|^{p-2}|z|+|z|^{p-1})\lesssim V_b^2(z).
\end{equation}

\subsubsection{The shifted cost}
Set $c(z) = \frac 1 p|z|^p$, so that $\nabla c(z)=|z|^{p-2}z$ and $c^\ast(y)=\frac1{p'}|y|^{p'}$. For $b\in\R^d$ we introduce
\begin{align}\label{eq:defhatc}
c_b(z) = c(z-b)-c(b)-\langle \nabla c(-b), z\rangle ,
\end{align}
which is the Bregman divergence of $c$ between $z-b$ and $-b$. In particular $c_b\ge0$, $c_b(0)=0$, $\nabla c_b(0)=0$, and
\begin{equation}\label{eq:gradhatc}
\nabla c_b(z)=\nabla c(z-b)-\nabla c(-b)=|z-b|^{p-2}(z-b)+|b|^{p-2}b.
\end{equation}

We record the classical estimate for Bregman divergences of powers: for every $q>1$ and $v,w\in\R^d$,
\begin{equation}\label{eq:bregmanpower}
\frac1q|v|^q-\frac1q|w|^q-|w|^{q-2}\langle w,v-w\rangle\simeq (|v|+|w|)^{q-2}|v-w|^2,
\end{equation}
with constants depending only on $d$ and $q$. Applying \eqref{eq:bregmanpower} with $q=p$, $v=u-b$, $w=v'-b$ gives the \emph{shifted convexity estimate}
\begin{equation}\label{eq:convexity}
c_b(u)-c_b(v)-\langle\nabla c_b(v),u-v\rangle\simeq W(u-b,v-b),\qquad W(x,y):=(|x|+|y|)^{p-2}|x-y|^2 .
\end{equation}

Taking $v=0$ in \eqref{eq:convexity} and using $|u-b|+|b|\simeq |u|+|b|$ we obtain
\begin{align}\label{eq:energyBoundV}
c_b(z) \simeq V_b^2(z).
\end{align}
By \eqref{eq:gradhatc} and the mean value theorem,
\begin{align}\label{eq:nablahatc}
|\nabla c_b(z)|\lesssim |b|^{p-2}|z|+|z|^{p-1},\qquad |\nabla c_b(x)-\nabla c_b(y)|\lesssim (|b|+|x|+|y|)^{p-2}|x-y|,
\end{align}
and consequently
\begin{align}\label{eq:continuity}
|c_b(x)-c_b(y)|\lesssim (|b|^{p-2}(|x|+|y|)+|x|^{p-1}+|y|^{p-1})|x-y|.
\end{align}
Since $\langle\nabla c_b(z),z\rangle\le|\nabla c_b(z)||z|\les V_b^2(z)\simeq c_b(z)$ and $V_b^2(\tau z)\simeq V_b^2(z)$ for $\tau\in[\frac12,2]$, we also have
\begin{equation}\label{eq:dilation}
|c_b(tz)-c_b(z)|\les |t-1|\,c_b(z)\qquad\text{for } t\in[\tfrac12,2].
\end{equation}

\subsubsection{The conjugate cost}
A direct computation gives
\begin{align}\label{eq:hatcstar}
c_b^\ast(\zeta) = \langle \zeta,b\rangle+c^\ast(\zeta+\nabla c(-b))-c^\ast(\nabla c(-b))
=\langle\zeta,b\rangle+\frac1{p'}\big|\zeta-|b|^{p-2}b\big|^{p'}-\frac1{p'}|b|^p,
\end{align}
and therefore
\begin{equation}\label{eq:gradhatcstar}
\nabla c_b^\ast(\zeta)=b+\big|\zeta-|b|^{p-2}b\big|^{p'-2}\big(\zeta-|b|^{p-2}b\big),\qquad \nabla c_b^\ast(0)=0 .
\end{equation}
Since $\nabla c^\ast(\nabla c(-b))=-b$, $c_b^\ast$ is the Bregman divergence of $c^\ast$ between $\zeta+\nabla c(-b)$ and $\nabla c(-b)$, so \eqref{eq:bregmanpower} with $q=p'$ gives
\begin{equation}\label{eq:energyBoundVstar}
c_b^\ast(\zeta)\simeq V_{b,\ast}^2(\zeta),\qquad
\langle \nabla c_b^\ast(\zeta_1)-\nabla c_b^\ast(\zeta_2),\zeta_1-\zeta_2\rangle\simeq \big(|\zeta_1-\eta_b|+|\zeta_2-\eta_b|\big)^{p'-2}|\zeta_1-\zeta_2|^2,
\end{equation}
where $\eta_b=|b|^{p-2}b$, and moreover
\begin{align}\label{lipnablacstar}
|\nabla c_b^\ast(\zeta_1)-\nabla c_b^\ast(\zeta_2)|\lesssim (|\zeta_1-\eta_b|+|\zeta_2-\eta_b|)^{p^\prime-2}|\zeta_1-\zeta_2|.
\end{align}

\subsubsection{A triangle inequality}
\begin{lemma}\label{triangle}
 For every $\eps\in(0,1)$ there exists $C_\eps>0$ such that for all $b,x,y\in\R^d$,
 \begin{equation}\label{eq:triangle}
  c_b(x+y)\le (1+\eps)c_b(x)+ C_\eps c_b(y).
 \end{equation}
\end{lemma}
\begin{proof}
For $b=0$ this is the classical inequality $|x+y|^p\le(1+\eps)|x|^p+C_\eps|y|^p$. For $b\neq0$ both sides are jointly homogeneous of degree $p$ in $(b,x,y)$, so we may assume $|b|=1$. By \eqref{eq:energyBoundV}, it is enough to prove that (possibly redefining the value of $\eps$)
 \begin{equation*}
  c_b(x+y)-c_b(x)\le \eps V_b^2(x)+ C_\eps V_b^2(y).
 \end{equation*}
By \eqref{eq:nablahatc} and Young's inequality we indeed obtain
\begin{multline*}
  c_b(x+y)-c_b(x)=\int_0^1 \langle\nabla c_b(x+ty), y\rangle\,dt\les  (|x|+|y| +|x|^{p-1}+|y|^{p-1})|y|\\
  \les \eps |x|^2+ \eps |x|^p + C_\eps(|y|^2+ |y|^p)\les \eps V_b^2(x)+ C_\eps V_b^2(y).
\end{multline*}
\end{proof}

\subsection{Optimal transport}\label{sec:prelimOT}

\subsubsection{Transport problems}
For finite measures $\lambda,\mu$ on $\R^d$ with $\lambda(\R^d)=\mu(\R^d)$, we write $\Pi(\lambda,\mu)$ for the set of couplings and
\begin{align*}
W_{c_b}(\lambda,\mu)= \inf_{\pi\in\Pi(\lambda,\mu)} \int c_b(y-x)\,d\pi,\qquad W_p^p(\lambda,\mu)=\inf_{\pi\in\Pi(\lambda,\mu)} \int |y-x|^p\,d\pi .
\end{align*}
Notice that since affine functions are null-Lagrangians, $\pi$ is optimal for $c_b$ if and only if 
it is optimal for the cost $(x,y)\mapsto c(y-x-b)$. For $\pi$ a transport plan we often write $u=y-x$ for the displacement.

\subsubsection{Change of coordinates}
The iteration we will apply shifts the source by a vector $a\in\R^d$. Hence, for a plan $\pi$ we let $\tilde\pi$ be the image of $\pi$ under $(x,y)\mapsto(\tilde x,\tilde y)=(x-a,y)$, with marginals $\tilde\lambda=\lambda(\cdot+a)$ and $\mu$.
 Then
\begin{equation}\label{eq:shiftrules}
\tilde u:=\tilde y-\tilde x=u+a,\qquad \tilde b:=b+a,\qquad \tilde u-\tilde b=u-b .
\end{equation}
By the discussion above we have:
\begin{lemma}\label{lem:shift}
If $\pi$ is optimal for $c_b$, then $\tilde\pi$ is optimal for $c_{b+a}$. \end{lemma}

\subsubsection{Energy and data}
We set, for $b\in\R^d$ and $R>0$, the  energy
\begin{align}\label{eq:defE}
E_{b,R}(\pi) = \frac 1 {R^d} \int_{\#_R} V_b^2(y-x)\,d\pi.
\end{align}
For a  measure $\nu$ and $R>0$, we define the $b$-free data term on balls centred at the origin by
\begin{align}\label{eq:defDat}
\Dat_R(\nu)=\frac1{R^d}W_p^p\big(\nu\rst B_R,\kappa_{\nu,B_R}\,dx\rst B_R\big)+R^p\,|\kappa_{\nu,B_R}-1|^p,
\end{align}
and $\Dat^\square_R(\nu)$ denotes the same quantity with $B_R$ replaced by the centred cube $Q_R$.  For a pair of measures we write $\Dat_R=\Dat_R(\lambda)+\Dat_R(\mu)$ and
\begin{equation}\label{eq:defD}
D_{b,R}=\Dat_R+|b|^{p-2}\,\Dat_R^{2/p}.
\end{equation}
The natural smallness condition which we will always assume to hold is
\begin{equation}\label{eq:smallnessphys}
R^{-p}\big(E_{b,R}+D_{b,R}\big)+R^{-1}|b|\ll1.
\end{equation}
For $p=2$ and $b=0$, $E_{0,R}$ and $\Dat_R$ are the energy and data terms of the harmonic approximation of \cite{GO} (see also \cite{goldman2021quantitative} and the lecture notes \cite{KoOt}); their analogues for $p$-costs are considered in \cite{koch2024geometric}.
The weighted term $|b|^{p-2}\Dat_R^{2/p}$ in $D_{b,R}$ reflects the two regimes of $V_b^2$, quadratic for $|z|\ll|b|$ and $p$-homogeneous for $|z|\gg|b|$, and is the price for working with the shifted cost $c_b$.
We record a simple comparison result between $W_{c_b}$ and $W_p$ whose proof is left to the reader.
\begin{lemma}\label{lem:weighteddata}
Let $\nu,\nu'$ be finite measures of mass $m$ and let $\gamma$ be optimal for $W_p(\nu,\nu')$. Then
\begin{equation}\label{eq:holderdata}
W_{c_b}(\nu,\nu')\le \int c_b(y-x)\,d\gamma\les W_p^p(\nu,\nu')+|b|^{p-2}\,m^{1-\frac2p}\,W_p^2(\nu,\nu').
\end{equation}
In particular, if $\nu=\lambda\rst B_R$, $\nu'=\kappa_{\lambda,B_R}dx\rst B_R$ and $\kappa_{\lambda,B_R}\le 2$, then
\[
R^{-d}W_{c_b}(\nu,\nu')\les \Dat_R(\lambda)+|b|^{p-2}\Dat_R(\lambda)^{2/p},
\qquad 
\textrm{and} \qquad  R^2|\kappa_{\lambda,B_R}-1|^2\le\Dat_R(\lambda)^{2/p}.\]
\end{lemma}

We now prove  triangle inequality and subadditivity for $W_{c_b}$. Recall that $c_b$ is not even so $W_{c_b}$ is not symmetric.  
\begin{lemma}\label{lem:triangleW}
For every $\delta\in(0,1)$ there exists $C_\delta>0$ such that for any finite measures $\lambda_1,\lambda_2,\lambda_3$ of equal mass,
 \begin{equation}\label{eq:triangleW}
 W_{c_b}(\lambda_1,\lambda_3)\leq \min\Big\{(1+\delta)W_{c_b}(\lambda_1,\lambda_2)+C_\delta W_{c_b}(\lambda_2,\lambda_3),\ C_\delta W_{c_b}(\lambda_1,\lambda_2)+(1+\delta) W_{c_b}(\lambda_2,\lambda_3)\Big\}.
 \end{equation}
 Moreover $W_{c_b}$ is subadditive: $W_{c_b}(\lambda_1+\lambda_1',\lambda_2+\lambda_2')\le W_{c_b}(\lambda_1,\lambda_2)+W_{c_b}(\lambda_1',\lambda_2')$.
\end{lemma}
\begin{proof}
Use the gluing lemma, see \cite[Chapter 1]{villani2009optimal} to glue optimal plans and apply Lemma \ref{triangle} pointwise. Subadditivity follows by adding plans.
\end{proof}

\subsubsection{Mass bounds from the data}
The following lemma converts smallness of the transport distance to a constant density into lower and upper bounds on the mass of balls,
 in the spirit of the arguments leading to the $L^\infty$ bounds in \cite{GO,goldman2021quantitative}.
\begin{lemma}\label{lem:mass}
There is $C_m\ge1$ such that the following holds. Let $\nu$ be a measure, $R>0$, $\kappa=\kappa_{\nu,R}\in[\frac12,2]$ (recall the definition \eqref{eq:kappa}) and $\mathcal W=W_p^p(\nu\rst B_R,\kappa\,dx\rst B_R)$. Let $B(z,r)\subset B_R$.
\begin{enumerate}
\item If $r^{d+p}\ge C_m\mathcal W$, then $\nu(B(z,r))\ge r^d/C_m$.
\item If $R\le 5$, $r\le1$, $B(z,2r)\subset B_R$ and $r^{p(d+1)}\ge \mathcal W$, then $\nu(B(z,r))\le C_mr^d$.
\end{enumerate}
\end{lemma}
\begin{proof}
Let $\gamma$ be optimal for $\mathcal W$ and let $\eta$ be a cutoff supported in $B(z,r)$, with $\eta=1$ on $B(z,r/2)$, $0\le\eta\le1$ and ${\rm Lip}\,\eta\le 2/r$. Then
\begin{align*}
\kappa|B_{r/2}|\le\int\eta\,\kappa\,dx&=\int\eta\,d\nu+\int(\eta(x')-\eta(x))\,d\gamma(x,x')\\
&\le \nu(B(z,r))+\frac2r\big(\nu(B(z,r))+\kappa|B_r|\big)^{1/p'}\mathcal W^{1/p},
\end{align*}
since only pairs with $x\in B(z,r)$ or $x'\in B(z,r)$ contribute. If $\nu(B(z,r))<\kappa|B_{r/2}|/2$, this gives $r^{d}\les r^{d/p'-1}\mathcal W^{1/p}$, that is $r^{d+p}\les\mathcal W$; this proves (1). For (2) use a cutoff equal to $1$ on $B(z,r)$, supported in $B(z,2r)$, with ${\rm Lip}\le1/r$. As above,
\[
\nu(B(z,r))\le\kappa|B_{2r}|+r^{-1}\big(\nu(B_R)+\kappa|B_R|\big)^{1/p'}\mathcal W^{1/p}\les r^d+r^{-1}\mathcal W^{1/p}\les r^d,
\]
where we used $\nu(B_R)=\kappa|B_R|\le 2|B_5|$ and $\mathcal W^{1/p}\le r^{d+1}$.
\end{proof}

\subsubsection{A Benamou--Brenier competitor}
 We now adapt a Benamou-Brenier type estimate from \cite[Lemma 7.4]{koch2024geometric} to our context.
\begin{lemma}\label{lem:BenamouBrenier}
Let $R>0$, let $\kappa_1,\kappa_2\in[\frac12,2]$ be constants and let $f,g$ be finite measures on $\partial B_R$ with $\kappa_1 |B_R|+f(\partial B_R)=\kappa_2 |B_R|+g(\partial B_R)$. Let $\xi\in L^1(B_R;\R^d)$ satisfy
\begin{align}\label{eq:BBeq}
\nabla \cdot \xi=\kappa_1-\kappa_2 \text{ in } B_R,\qquad \xi\cdot \nu = g-f \text{ on } \partial B_R,
\end{align}
in the weak sense. Let $s_t= (1-t) \kappa_1+t\kappa_2$. Then
\begin{align*}
W_{c_b}(\kappa_1 d x\llcorner B_R+f,\kappa_2 d x\llcorner B_R+g)&\leq \int_{B_R}\int_0^1 s_t\,c_b(s_t^{-1}\xi)\,dt\,dx\\
&\le \int_{B_R} c_b(\xi)\,dx+C\,\max_i|\kappa_i-1|\int_{B_R} c_b(\xi)\,dx .
\end{align*}
\end{lemma}
\begin{proof}
Set $\rho_t=s_t\,dx\rst B_R+(1-t)f+tg$ and $j_t=\xi\,dx\rst B_R$. For $\zeta\in C^1_c(\R^d)$,
\[
\frac{d}{dt}\int\zeta\,d\rho_t=(\kappa_2-\kappa_1)\int_{B_R}\zeta+\int\zeta\,d(g-f)=\int_{B_R}\langle\xi,\nabla\zeta\rangle ,
\]
by \eqref{eq:BBeq}; hence $\partial_t\rho_t+\nabla\cdot j_t=0$ in $\R^d$, with $\rho_0=\kappa_1dx\rst B_R+f$ and
 $\rho_1=\kappa_2dx\rst B_R+g$.  The velocity $v_t=j_t/\rho_t$ equals $\xi/s_t$ in $B_R$ and $0$ on $\partial B_R$, and $\int_0^1\!\int|v_t|\,d\rho_t\,dt<\infty$. By the superposition principle \cite[Theorem 8.2.1]{AGS}, $\rho_t$ is the time marginal of a probability measure on absolutely continuous curves solving $\dot\gamma=v_t(\gamma)$; coupling $\gamma(0)$ with $\gamma(1)$ and using Jensen's inequality $c_b(\gamma(1)-\gamma(0))\le\int_0^1c_b(\dot\gamma)\,dt$ gives the first inequality. The second follows from $sc_b(\xi/s)=c_b(\xi)+(s-1)c_b(\xi)+s(c_b(\xi/s)-c_b(\xi))$ and \eqref{eq:dilation}.
\end{proof}

\subsubsection{Localization}
The following localization result is essentially contained in \cite[Lemma 6.1]{koch2024geometric}, see also \cite{goldman2021quantitative}.

\begin{lemma}\label{inp:loc}
There are $\eps_{\rm loc}>0$ and $C_{\rm loc}<\infty$ such that the following holds. Let $\nu$ be a measure,  $U$ be either the ball $B_R$ or the cube $Q_R$, and let 
$z\in\R^d$, $\rho\ge R/16$ with $B_{2\rho}(z)\subset U$. If $R^{-p}\Dat_U(\nu)\le\eps_{\rm loc}$, then
\[
\frac1\rho\int_\rho^{2\rho}\Dat_r(\nu(\cdot+z))\,dr\le C_{\rm loc}\,\Dat_U(\nu).
\]
\end{lemma}

\section{The deterministic input}\label{sec:deterministic}

In this section all statements are deterministic.  Throughout, $\lambda$ and $\mu$ are finite measures on $\R^d$ with compact support and equal mass, and $\pi\in\Pi(\lambda,\mu)$ is optimal for $c_b$. We remind the reader that $c_b$ was defined in Section \ref{sec:prelimOT}. At unit scale we fix $R'\in[4,5]$ and write
\[
E=E_{b,4}(\pi),\qquad \Dat=\Dat_{R'}(\lambda)+\Dat_{R'}(\mu),\qquad D=D_{b,R'}=\Dat+|b|^{p-2}\Dat^{2/p},\qquad \hE=E+D .
\]
Recall the definition \eqref{def:OmegaR} of $\Omega_R$. The data radius $R'$ is allowed to vary in $[4,5]$ since due to e.g. Lemma \ref{inp:loc} good estimates on $D$ only hold for some radii; all constants are uniform in $R'\in[4,5]$. If $\Dat\le 1$ then $|\kappa_{\lambda,R'}-1|,|\kappa_{\mu,R'}-1|\le\frac14\Dat^{1/p}\le\frac14$.

\subsection{Harmonic approximation}\label{sec:harm}

\subsubsection{A local $L^\infty$ bound}
For $R\in(2,3)$ we let $\mathcal C_R=\{(x,y)\in\Omega_R\colon x\notin B_R\text{ or }y\notin B_R\}$ be the set of pairs whose trajectory meets $\overline{B_R}$ without being contained in $B_R$.
The following lemma is an adaptation to our setting of $L^\infty$ bounds on the displacement in optimal transport, see \cite{bouchitte2007new, GO, goldman2021quantitative, koch2024geometric,goldman2025partial}.
\begin{lemma}\label{lem:Linfty}
There exist $\eps_\infty>0$ and $C_\infty<\infty$ such that, if $|b|+E+\Dat\le\eps_\infty$, then:
\begin{enumerate}
\item every $(x,y)\in\spt\pi$ whose trajectory meets $\overline{B_3}$ satisfies $x\in B_{7/2}$ or $y\in B_{7/2}$;
\item for every $(x,y)\in\spt\pi\cap\#_{7/2}$,
\begin{equation}\label{eq:LinftyBound}
|y-x|\le M:=C_\infty\big(E+\Dat+|b|^p\big)^{\frac1{d+p}}+|b| ;
\end{equation}
\item with $\mathcal C_R$ as above,
\begin{equation}\label{eq:crossing}
\int_2^3\pi(\mathcal C_R)\,dR\les M,\qquad \int_2^3\int_{\mathcal C_R}c_b(y-x)\,d\pi\,dR\les ME .
\end{equation}
\end{enumerate}
\end{lemma}

The proof  is an adaptation to the shifted cost of the  estimate \cite[Lemma 2.3]{d2026concentration}, 
combined with the counting argument of \cite[Lemma 2.4]{huesmann2024there}. 

\begin{proof}

We use the original displacement $w=y-x-b$ (see \eqref{eq:shiftrules}). By the cyclical monotonicity of $\spt\pi$, for $(x,y),(x',y')\in\spt\pi$,
\begin{equation}\label{eq:twoopt}
|y-x-b|^p+|y'-x'-b|^p\le |y'-x-b|^p+|y-x'-b|^p .
\end{equation}
Since $V_b^2(u)\ge|u|^p$ and $\pi(\#_4)\le\lambda(B_4)+\mu(B_4)\les |B_5|$,
\begin{equation}\label{eq:wenergy}
\int_{\#_4}|y-x-b|^p\,d\pi\les\int_{\#_4}(|y-x|^p+|b|^p)\,d\pi\les E+|b|^p.
\end{equation}
Finally, $W_p^p(\lambda\rst B_{R'},\kappa_{\lambda,R'}dx\rst B_{R'})\les\Dat$, and similarly for $\mu$.

\emph{Proof of (2).} Let $(x,y)\in\spt\pi$ with $x\in B_{7/2}$ (the case $y\in B_{7/2}$ is similar), put $\ell=|w|$ and $e=w/\ell$ (if $\ell=0$ there is nothing to prove). 
Let $t_0=(E+\Dat+|b|^p)^{1/(d+p)}$, so that the claim follows from
\begin{equation}\label{claimLinfty2}
  \ell \les t_0.
\end{equation}
Let $K\ge1$ be a constant fixed below and set $s'=Kt_0$. For $\eps_\infty$ small enough depending on $K$, $s'\le\frac18$. Let $C_2\ge2$ be a constant fixed below.
We may assume that  $\ell\ge C_2s'$ since otherwise there is nothing to prove. Consider $B'=B(x+s'e,s'/2)\subset B_{4}$ and $(x',y')\in\spt\pi$ with
 $x'\in B'$, and put $\ell'=|y'-x'-b|$. 
 
 We claim that $\ell'>s'$. Indeed, if $\ell'\le s'$, then
 writing $x'-x=s'e+e'$ with $|e'|\le s'/2$, we have $|y-x'-b|=|w-s'e-e'|\le\ell-\frac{s'}2$ and $|y'-x-b|\le \ell'+\frac32s'$. Since we assumed that  
 $\ell'\le s'$, \eqref{eq:twoopt} gives
\[
\frac{p s'}2\Big(\frac\ell2\Big)^{p-1}\le \ell^p-\Big(\ell-\frac{s'}2\Big)^p\le \Big(\ell'+\frac32s'\Big)^p-\ell'^p\le\Big(\frac52 s'\Big)^p,
\]
hence $\ell\le C_2's'$ with $C_2'=C_2'(p)$. We choose $C_2=2C_2'+2$ which gives a contradiction to the initial assumption $\ell\ge C_2s'$ unless $\ell'>s'$. 
Thus, every $(x',y')\in\spt\pi$ with $x'\in B'$ satisfies $\ell'>s'$, and \eqref{eq:wenergy} yields
\[
\lambda(B')\le (s')^{-p}\int_{\#_4}|y'-x'-b|^p\,d\pi\les (s')^{-p}t_0^{d+p}.
\]
On the other hand, Lemma \ref{lem:mass}(1) with $r=s'/2$ gives $\lambda(B')\ge (s'/2)^d/C_m$ provided $(s'/2)^{d+p}\ge C_m5^dt_0^{d+p}$, which holds if $K$ is
 large enough. Both bounds are incompatible if $K^{d+p}$ is large enough. We fix such a $K$ and conclude $\ell<C_2Kt_0$, so $|y-x|\le\ell+|b|\le M$. 

\emph{Proof of (1).} Let $(x,y)\in\spt\pi$ with $|x|,|y|\ge\frac72$ and $|X(t)|\le3$ for some $t\in(0,1)$. Then $t|y-x|\ge\frac12$ and $(1-t)|y-x|\ge\frac12$, so that, with $w,\ell$ as above and $|b|\le\frac18$, $a:=t\ell\ge\frac14$ and $a':=(1-t)\ell\ge\frac14$. Let $\tilde X=x+tw$, so $|\tilde X-X(t)|\le|b|$, and let $\rho\in(0,\frac14]$ be a constant fixed below. For $(x',y')\in\spt\pi$ with $x'\in B(\tilde X,\rho)\subset B_{7/2}$ and $\ell'=|y'-x'-b|\le\rho$, \eqref{eq:twoopt} gives
\[
(a+a')^p\le (a'+\rho)^p+(a+2\rho)^p .
\]
Since $p\ge2$, $(a+a')^p-a^p-a'^p\ge c_p\,aa'(a+a')^{p-2}\ge c\,\ell^{p-1}$, while $(a'+\rho)^p-a'^p+(a+2\rho)^p-a^p\le C\rho\,\ell^{p-1}$. For $\rho$ small this is a contradiction. Hence all such $(x',y')$ have $\ell'>\rho$, and as above $\rho^d/C_m\le\lambda(B(\tilde X,\rho))\les\rho^{-p}(E+|b|^p)$ once $\rho^{d+p}\ge C_m5^d\Dat$. This is impossible if $\eps_\infty$ is small (depending on the fixed $\rho$).

\emph{Proof of (3).} By (1) and (2), every pair in $\Omega_R$, $R\le3$, lies in $\#_{7/2}$ and satisfies $|y-x|\le M$. If $(x,y)\in\mathcal C_R$ then $\min_t|X(t)|\le R\le\max(|x|,|y|)$, so the set of $R$ with $(x,y)\in\mathcal C_R$ is contained in an interval of length at most $|y-x|\le M$. By Fubini's theorem, $\int_2^3\pi(\mathcal C_R)dR\le M\pi(\#_4)$ and $\int_2^3\int_{\mathcal C_R}c_b\,d\pi\,dR\le M\int_{\#_4}c_b\,d\pi\les ME$ by \eqref{eq:energyBoundV}.
\end{proof}

\subsubsection{Boundary fluxes}
For $R\in(2,3)$ we define the incoming and outgoing fluxes $f_R,g_R$ on $\partial B_R$ by
\begin{align*}
\int_{\partial B_R} \zeta\, d f_R = \int_{\Omega_R\cap\{x\notin B_R\}} \zeta(X(t_-))\,d \pi,\qquad \int_{\partial B_R} \zeta \,d g_R = \int_{\Omega_R\cap\{y\notin B_R\}} \zeta(X(t_+))\,d\pi .
\end{align*}
If $\lambda(\partial B_R)=\mu(\partial B_R)=0$, the image of $\pi\rst\Omega_R$ under $(x,y)\mapsto(X(t_-),X(t_+))$ has marginals $\lambda\rst B_R+f_R$ and $\mu\rst B_R+g_R$, and in particular
\begin{equation}\label{eq:massbalance}
\lambda(B_R)+f_R(\partial B_R)=\mu(B_R)+g_R(\partial B_R).
\end{equation}
The following lemma is a simple extension of the projection argument from \cite{goldman2021quantitative}, see also \cite[Lemma 2.3]{koch2024geometric}. We do not repeat the proof.
\begin{lemma}\label{lem:bathtub}
Let $\Lambda>0$, $R\in[2,3]$, let $g'$ be a measure on $B_5$ with $g'\le\Lambda\,dx$, and let $\bar g=(\Pi_R)_\#g'$ with $\Pi_R(z)=Rz/|z|$. Then $\bar g$ has a density 
with respect to the surface measure on $\partial B_R$ and, for every $k\ge1$,
\[
\int_{\partial B_R}\bar{g}^k\le C_k\Lambda^{k-1}\int\big|R-|z|\big|^{k-1}\,dg'(z).
\]
\end{lemma}
The following result is an adaptation of \cite[Lemma 5.1]{koch2024geometric}, see also \cite{goldman2021quantitative,KoOt}, to the shifted cost $c_b$.

\begin{lemma}\label{lem:boundaryData}
Assume $|b|+E+\Dat\le\eps_\infty$. For a.e. $R \in (2,3)$, there exist measures $\bar f_R,\bar g_R$ on $\partial B_R$,  with $\bar f_R(\partial B_R)=f_R(\partial B_R)$ and $\bar g_R(\partial B_R)=g_R(\partial B_R)$, having densities with respect to the surface measure, such that
\begin{align}\label{eq:bdryData}
\int_2^3 \Big(W_{c_b}(f_R,\bar f_R)+ W_{c_b}(g_R,\bar g_R)\Big)dR\lesssim M \hE,\qquad
\int_2^3 \int_{\partial B_R} \big(V_b^2(\bar f_R)+V_b^2(\bar g_R)\big)\, d R \lesssim \hE,
\end{align}
with $\hE$ defined in \eqref{eq:defE} and $M$ in \eqref{eq:LinftyBound}.
\end{lemma}
\begin{proof}
We construct $\bar g_R$; the construction of $\bar f_R$ is symmetric. Let $\gamma$ be optimal for $W_p(\mu\rst B_{R'},\kappa_{\mu,R'}dz\rst B_{R'})$ and disintegrate $\gamma(dy,dz)=\mu\rst B_{R'}(dy)\,\gamma(dz|y)$. By Lemma \ref{lem:weighteddata},
\begin{equation}\label{eq:dataplan}
\int V_b^2(z-y)\,d\gamma\les D .
\end{equation}
For pairs in $\Omega_R\cap\{y\notin B_R\}$ we have $y\in B_{R+M}\subset B_{R'}$ by Lemma \ref{lem:Linfty}. Let $\tilde\pi_R(dx,dy,dz)=\pi\rst(\Omega_R\cap\{y\notin B_R\})(dx,dy)\,\gamma(dz|y)$, let $g'_R$ be its $z$-marginal and $\bar g_R=(\Pi_R)_\#g'_R$. Since the $y$-marginal of $\pi\rst(\Omega_R\cap\{y\notin B_R\})$ is at most $\mu\rst B_{R'}$, we have $g'_R\le\kappa_{\mu,R'}\,dz\le 2\,dz$.

The image of $\tilde\pi_R$ under $(x,y,z)\mapsto(X(t_+),\Pi_Rz)$ is a coupling of $g_R$ and $\bar g_R$. Since $|X(t_+)|=R$, we have $|\Pi_Rz-z|=||z|-R|\le|z-X(t_+)|$, hence $|\Pi_Rz-X(t_+)|\le 2|z-X(t_+)|\le2(|z-y|+|y-x|)$ and, by \eqref{eq:energyBoundV}, \eqref{eq:triangleV} and \eqref{eq:scalingV},
\[
W_{c_b}(g_R,\bar g_R)\les\int\big(V_b^2(y-x)+V_b^2(z-y)\big)\,d\tilde\pi_R .
\]
For fixed $(x,y)$, the set of $R\in(2,3)$ with $(x,y)\in\Omega_R$ and $y\notin B_R$ is contained in $[\min_t|X(t)|,|y|]$, an interval of length at most $|y-x|\le M$. By Fubini's theorem, \eqref{eq:dataplan} and $\int_{\#_4}V_b^2(y-x)d\pi\les E$ we obtain the first bound in \eqref{eq:bdryData}.

For the second bound we use Lemma \ref{lem:bathtub} with $\Lambda=2$ and $k\in\{2,p\}$, together with $V_b^2(h)\les|b|^{p-2}h^2+h^p$:
\begin{align*}
\int_{\partial B_R}V_b^2(\bar g_R)&\les\int\big(|b|^{p-2}\big|R-|z|\big|+\big|R-|z|\big|^{p-1}\big)dg_R'(z)\\
&\les\int\big(|b|^{p-2}\ell+\ell^{p-1}\big)\,d\tilde\pi_R,\qquad \ell=|z-y|+|y-x| .
\end{align*}
Integrating in $R$ and using again that the relevant set of $R$ has length at most $|y-x|$, together with $|y-x|\ell^{k-1}\les |y-x|^k+|z-y|^k$, we get
\[
\int_2^3\int_{\partial B_R}V_b^2(\bar g_R)\,dR\les\int\Big(|b|^{p-2}(|y-x|^2+|z-y|^2)+|y-x|^p+|z-y|^p\Big)d\tilde\pi\les E+D,
\]
where $\tilde\pi=\pi\rst\#_4(dx,dy)\gamma(dz|y)$ and we used \eqref{eq:dataplan}.
\end{proof}

\subsubsection{Choice of a good radius}\label{sec:goodRadius}
Assume $|b|+E+\Dat\le\eps_\infty$ and $\Dat\le\eps_{\rm loc}$. By Lemma \ref{lem:Linfty}(3), Lemma \ref{lem:boundaryData}, Lemma \ref{inp:loc} (applied to $\lambda$ and $\mu$ with $U=B_{R'}$, $z=0$, $\rho=2$) and Chebyshev's inequality, there is $R\in(2,3)$ such that $\lambda(\partial B_R)=\mu(\partial B_R)=0$ and
\begin{align}
\pi(\mathcal C_R)&\lesssim M,\label{eq:measureCrossing}\\
\int_{\mathcal C_R} c_b(y-x)\,d\pi&\lesssim M E,\label{eq:crossingTrajectoriesCost}\\
W_{c_b}(f_R,\bar f_R)+W_{c_b}(g_R,\bar g_R)&\lesssim M \hE,\label{eq:WassersteinData}\\
\int_{\partial B_R} \big(V_b^2(\bar f_R)+V_b^2(\bar g_R)\big)&\lesssim \hE,\label{eq:IntegrabilityData}\\
\Dat_R(\lambda)+\Dat_R(\mu)\les\Dat,\qquad D_{b,R}&\les D .\label{eq:goodData}
\end{align}

\subsubsection{The harmonic approximation theorem}
We use the results of Appendix \ref{app:PDE} with $B=B_R$: for a density $h$ on $\partial B_R$ with $\int V_b^2(h)<\infty$ and $\omega=|B_R|^{-1}\int_{\partial B_R} h$, $\phi$ denotes the weak solution of
\begin{equation}\label{eq:PDEmain}
\nabla\cdot\nabla c_b^\ast(\nabla\phi)=\omega\ \text{ in }B_R,\qquad \nabla c_b^\ast(\nabla\phi)\cdot\nu=h\ \text{ on }\partial B_R ,
\end{equation}
$\xi=\nabla c_b^\ast(\nabla\phi)$, and $\phi_r,\xi_r$ denote the corresponding objects with $h$ replaced by $P_rh$, the average of $h$ over rotations of size at most $r$ (Definition \ref{def:Pr}).

\begin{theorem}\label{harmapprox}
Let $\beta\in(0,1)$ be the exponent of Theorem \ref{thm:regularity}(iii)--(iv). There is $C_0<\infty$, and for every $\delta>0$ there exist $\eps>0$ and $C_\delta<\infty$, such that the following holds. If $\pi$ is optimal for $c_b$, $R'\in[4,5]$ and
 \begin{equation}\label{hypsmall}
  |b|+E_{b,4}(\pi)+ D_{b,R'}\le \eps ,
 \end{equation}
then there exist $R\in(2,3)$, $\omega\in\R$ and $\xi=\nabla c_b^*(\nabla \phi)$ with $\nabla\cdot \xi=\omega$ in $B_R$, $\xi\in C^{0,\beta}(\overline{B_R})$,
\begin{equation}\label{eq:estimxithm}
 V_b^2\big(\|\xi\|_{C^{0,\beta}(B_{7/4})}\big)\le C_0\big( E_{b,4}(\pi)+ D_{b,R'}\big),
\end{equation}
and recalling the definition \eqref{eq:convexity} of $W$,
\begin{equation}\label{conclharmapprox}
\begin{aligned}
 \int_{\#_1} W(y-x-b,\xi(x)-b)\,d\pi(x,y)&
 \le \delta E_{b,4}(\pi) + C_\delta D_{b,R'}.
\end{aligned}
\end{equation}
\end{theorem}

For $(x,y)\in\#_1\cap\spt\pi$ we have $x\in B_{1+M}\subset B_{7/4}$ by Lemma \ref{lem:Linfty}, so $\xi(x)$ in \eqref{conclharmapprox} is well defined.

\begin{proof}[Proof of Theorem \ref{harmapprox}]
Fix $\delta\in(0,1)$; the parameter $r\in(0,\frac1{12}]$ will be chosen depending on $\delta$, and then $\eps$ depending on $\delta$ and $r$. We assume $\eps\le\min(\eps_\infty,\eps_{\rm loc})$, choose $R$ as in Section \ref{sec:goodRadius}, and write $\Omega=\Omega_R$, $f=f_R$, $g=g_R$, $\bar f=\bar f_R$, $\bar g=\bar g_R$, $\kappa_\lambda=\kappa_{\lambda,R}$, $\kappa_\mu=\kappa_{\mu,R}$. Let
\[
h=\bar g-\bar f,\qquad \omega=\frac1{|B_R|}\int_{\partial B_R}h=\kappa_\lambda-\kappa_\mu ,
\]
where the last identity follows from \eqref{eq:massbalance} and $\bar f(\partial B_R)=f(\partial B_R)$, $\bar g(\partial B_R)=g(\partial B_R)$. Let $\phi,\xi,\phi_r,\xi_r$ be as above. By \eqref{eq:IntegrabilityData} and \eqref{eq:triangleV}, $\int_{\partial B_R}V_b^2(h)\les\hE$. Theorem \ref{thm:regularity} then gives
\begin{equation}\label{eq:regapplied}
\int_{B_R}c_b(\xi)\les\hE,\qquad V_b^2\big(\|\xi_r\|_{C^{0,\beta}(B_{7/4})}\big)\le C\hE,\qquad V_b^2\big(\|\xi_r\|_{C^{0,\beta}(\overline{B_R})}\big)\le C_r\hE.
\end{equation}
We emphasize the fact that  the second constant does not depend on $r$. We take $\xi_r$ as the field in the statement; \eqref{eq:estimxithm} is the second estimate in \eqref{eq:regapplied}. Since $\nabla\phi_r=\nabla c_b(\xi_r)$, \eqref{eq:nablahatc} and \eqref{eq:elementary} give $\phi_r\in C^1(\overline{B_R})$ and
\begin{equation}\label{eq:gradphir}
\sup_{\overline{B_R}}V_{b,\ast}^2(\nabla\phi_r)\les \sup_{\overline{B_R}}V_b^2(\xi_r)\le C_r\hE .
\end{equation}

\step 1 {The main identity.} By \eqref{eq:convexity} with $u=\dot X$ and $v=\xi_r(X)$, and since $\nabla c_b(\xi_r)=\nabla\phi_r$,
\begin{align*}
\int_\Omega \int_{t_-}^{t_+} W(\dot X-b,\xi_r(X)-b)\, d t\, d\pi\lesssim{}& \int_\Omega c_b(y-x)\,d\pi-\int_\Omega \int_{t_-}^{t_+} c_b(\xi_r(X))\, d t\, d\pi\\
&-\int_\Omega \int_{t_-}^{t_+} \langle \nabla \phi_r(X),\dot X-\xi_r(X)\rangle\, d t\, d\pi ,
\end{align*}
where we used $\int_{t_-}^{t_+}c_b(\dot X)dt\le c_b(y-x)$. Since $\frac d{dt}\phi_r(X)=\langle\nabla\phi_r(X),\dot X\rangle$ and $(X(t_-),X(t_+))$ has marginals $\lambda\rst B_R+f$ and $\mu\rst B_R+g$,
\begin{align*}
\int_\Omega \int_{t_-}^{t_+} \langle\nabla \phi_r(X),\dot X\rangle\,dt\,d\pi&= \int_\Omega \big(\phi_r(X(t_+))-\phi_r(X(t_-))\big)d \pi\\
&= \int_{\partial B_R} \phi_r\,d(g- f)+\int_{B_R} \phi_r\, d(\mu-\lambda).
\end{align*}
Writing $F=c_b(\xi_r)-\langle\nabla c_b(\xi_r),\xi_r\rangle$ and testing the equation of $\phi_r$ with $\phi_r$,
\begin{align*}
\int_{B_R} \langle \nabla \phi_r,\xi_r\rangle\, d x= \int_{\partial B_R} \phi_r\,P_rh-\omega\int_{B_R} \phi_r \,dx.
\end{align*}
Combining, we arrive at
\begin{align}\label{eq:estimateMain}
&\int_\Omega \int_{t_-}^{t_+} W(\dot X-b,\xi_r(X)-b)\, d t\, d\pi\nonumber\\
&\qquad\lesssim \Big[\int_\Omega c_b(y-x)\,d\pi-\int_{B_R} c_b(\xi_r)\,d x\Big]
- \Big[\int_\Omega \int_{t_-}^{t_+} F(X)\,d t\, d\pi-\int_{B_R} F\, d x\Big]\nonumber\\
&\qquad\quad+ \int_{\partial B_R} \phi_r\,\big[P_rh-(g-f)\big] +\Big[\int_{B_R} \phi_r \,d(\lambda-\kappa_\lambda dx)-\int_{B_R}\phi_r\,d(\mu-\kappa_\mu dx)\Big]\nonumber\\
&\qquad=:T_1+T_2+T_3+T_4 ,
\end{align}
where we used $\omega=\kappa_\lambda-\kappa_\mu$ in the last bracket.


\step 2 {The first term.} By convexity of $c_b$ and $\nabla c_b(\xi)=\nabla\phi$,
\[
\int_{B_R}c_b(\xi)-\int_{B_R}c_b(\xi_r)\le\int_{B_R}\langle\nabla\phi,\xi-\xi_r\rangle=\int_{\partial B_R}\phi\,(h-P_rh)\les r^{\frac{p'-1}2}\hE ,
\]
where the identity uses that $\xi$ and $\xi_r$ have the same divergence $\omega$, and the inequality is Theorem \ref{thm:regularity}(ii). Next, $\pi\rst\Omega$ is optimal between its marginals $\lambda_\Omega,\mu_\Omega$. The plan $(X(t_-),x)_\#\pi\rst\Omega$ couples $\lambda\rst B_R+f$ with $\lambda_\Omega$ at cost at most $\int_{\mathcal C_R}c_b(t_-(y-x))d\pi\le\int_{\mathcal C_R}c_b(y-x)d\pi$ (convexity and $c_b(0)=0$), and similarly on the target side. By Lemma \ref{lem:triangleW} and \eqref{eq:crossingTrajectoriesCost},
\[
\int_\Omega c_b(y-x)\,d\pi\le (1+\delta) W_{c_b}(\lambda \llcorner B_R+f,\mu\llcorner B_R+g)+ C_\delta M E .
\]
Applying Lemma \ref{lem:triangleW} twice more, subadditivity, \eqref{eq:WassersteinData}, \eqref{eq:goodData} and Lemma \ref{lem:weighteddata},
\begin{align*}
W_{c_b}(\lambda \llcorner B_R+f,\mu\llcorner B_R+g)&\leq (1+\delta)^2 W_{c_b}(\kappa_{\lambda} d x \llcorner B_R+\bar f,\kappa_{\mu}d x \llcorner B_R+\bar g)+ C_\delta(M\hE+D).
\end{align*}
By Lemma \ref{lem:BenamouBrenier} (with $\kappa_1=\kappa_\lambda$, $\kappa_2=\kappa_\mu$, and $\xi$, which satisfies \eqref{eq:BBeq}), $|\kappa_\lambda-1|+|\kappa_\mu-1|\les D^{1/p}$ and \eqref{eq:regapplied},
\[
W_{c_b}(\kappa_{\lambda}dx \llcorner B_R+\bar f, \kappa_{\mu}d x \llcorner B_R+\bar g)\le(1+CD^{1/p})\int_{B_R}c_b(\xi)\le\int_{B_R}c_b(\xi)+CD^{1/p}\hE .
\]
Collecting,
\begin{equation}\label{eq:T1}
T_1\les\big(\delta+D^{1/p}+r^{\frac{p'-1}2}\big)\hE+C_\delta\big(M\hE+D\big).
\end{equation}

\step 3 {The second term.} Since $B_R$ is convex, $X(t)\in\overline{B_R}$ exactly for $t\in[t_-,t_+]$, so
\[
\int_\Omega\int_{t_-}^{t_+} F(X)\,dt\,d\pi-\int_{B_R}F\,d\lambda=\int_\Omega\int_0^1\big(\one_{\overline{B_R}}(X(t))F(X(t))-\one_{B_R}(x)F(x)\big)dt\,d\pi .
\]
For pairs outside $\mathcal C_R$ the integrand is bounded by $[F]_{C^{0,\beta}(B_R)}|y-x|^\beta$, and for pairs in $\mathcal C_R$ by $2\sup|F|$. Moreover, if $\gamma_\lambda$ is optimal for $W_p(\lambda\rst B_R,\kappa_\lambda dx\rst B_R)$, then by Jensen's inequality $\int|x-x'|^\beta d\gamma_\lambda\les \Dat_R(\lambda)^{\beta/p}$. Hence, by \eqref{eq:LinftyBound}, \eqref{eq:measureCrossing} and \eqref{eq:goodData},
\begin{align*}
|T_2|&\le 2\sup_{B_R}|F|\,\pi(\mathcal C_R)+[F]_{C^{0,\beta}(B_R)}M^\beta\pi(\Omega)+\Big|\int_{B_R}F\,d(\lambda-\kappa_\lambda dx)\Big|+|\kappa_\lambda-1|\int_{B_R}|F|\\
&\les \sup_{B_R}|F|\,(M+D^{1/p})+[F]_{C^{0,\beta}(B_R)}\big(M^\beta+D^{\beta/p}\big).
\end{align*}
By \eqref{eq:energyBoundV} and \eqref{eq:nablahatc}, $|F|\les V_b^2(\xi_r)$. With $\xi_x=\xi_r(x)$, \eqref{eq:continuity} and \eqref{eq:nablahatc} give for $x,y\in B_R$
\begin{align*}
|F(x)-F(y)|\les \big(|b|^{p-2}(|\xi_x|+|\xi_y|)+|\xi_x|^{p-1}+|\xi_y|^{p-1}\big)[\xi_r]_{C^{0,\beta}(B_R)}|x-y|^\beta,
\end{align*}
where we used $|\xi_y|^{p-2}|\xi_x|\le|\xi_x|^{p-1}+|\xi_y|^{p-1}$. By \eqref{eq:YoungsV} with $\eps=1$ and \eqref{eq:elementary}, $[F]_{C^{0,\beta}(B_R)}\les V_b^2([\xi_r]_{C^{0,\beta}})+\sup V_b^2(\xi_r)$. By \eqref{eq:regapplied},
\begin{equation}\label{eq:T2}
|T_2|\le C_r\big(M+M^\beta+D^{\beta/p}\big)\hE .
\end{equation}

\step 3 {The third term.} We split $T_3=\int\phi_r(P_rh-h)+\int\phi_r\,d(\bar g-g)-\int\phi_r\,d(\bar f-f)$. By Theorem \ref{thm:regularity}(ii) applied to $\phi_r$, the first term is $\les r^{(p'-1)/2}\hE$. Let $\gamma_g$ be optimal for $W_{c_b}(g,\bar g)$; its mass is $m_g=g(\partial B_R)\le\pi(\mathcal C_R)\les M\le1$. For $\eta\in(0,1]$, by \eqref{eq:YoungsV}, Jensen's inequality, $V_b^2(m t)\le m^2V_b^2(t)$ for $m\le1$, and \eqref{eq:WassersteinData},
\begin{align*}
\Big|\int\phi_r\,d(\bar g-g)\Big|&\le \sup_{\overline{B_R}}|\nabla\phi_r|\int|y'-y|\,d\gamma_g\\
&\le \eta\sup V_{b,\ast}^2(\nabla\phi_r)+C\eta^{1-p}V_b^2\Big(\int|y'-y|d\gamma_g\Big)\\
&\le \eta C_r\hE+C\eta^{1-p}m_g\int V_b^2(y'-y)\,d\gamma_g\le \eta C_r\hE+C\eta^{1-p}M^2\hE .
\end{align*}
The same holds for $f$. Choosing $\eta=\delta/C_r$,
\begin{equation}\label{eq:T3}
|T_3|\les\big(r^{\frac{p'-1}2}+\delta\big)\hE+C_{\delta,r}M^2\hE .
\end{equation}

\step 4 {The fourth term.} With $\gamma_\lambda$ optimal for $W_{c_b}(\lambda\llcorner B_R,\kappa_{\lambda} dx \llcorner B_R)$, $\int|x-x'|d\gamma_\lambda\les\Dat_R(\lambda)^{1/p}$, and $V_b^2(C\Dat_R(\lambda)^{1/p})\les D_{b,R}\les D$. As in Step 3',
\[
\Big|\int_{B_R}\phi_r\,d(\lambda-\kappa_\lambda dx)\Big|\le\sup|\nabla\phi_r|\int|x-x'|d\gamma_\lambda\le\eta C_r\hE+C\eta^{1-p}D ,
\]
and similarly for $\mu$. With $\eta=\delta/C_r$,
\begin{equation}\label{eq:T4}
|T_4|\le \delta\hE+C_{\delta,r}D .
\end{equation}

\step 5 {The left-hand side.} Let $(x,y)\in\#_1\cap\spt\pi$ and $u=y-x$. By Lemma \ref{lem:Linfty}, $|u|\le M$. Hence the segment $[x,y]$ lies in $B_{1+M}\subset B_{7/4}$, and $t_-=0$, $t_+=1$. Let $\Xi=\|\xi_r\|_{C^{0,\beta}(B_{7/4})}$, so that $V_b^2(\Xi)\les\hE$ by \eqref{eq:regapplied}. Fix $t\in[0,1]$, set $a=\xi_r(X(t))$, $a_0=\xi_r(x)$ and $\delta_t=a-a_0$, so $|a|\le\Xi$ and $|\delta_t|\le\Xi M^\beta$. Since 
\[(|u-b|+|a_0-b|)^{p-2}\le(|u-b|+|a-b|+|\delta_t|)^{p-2}\] 
and $|u-a_0|^2\le2|u-a|^2+2|\delta_t|^2$,
\[
W(u-b,a_0-b)\les W(u-b,a-b)+|\delta_t|^{p-2}|u-a|^2+(|u-b|+|a-b|+|\delta_t|)^{p-2}|\delta_t|^2 .
\]
By Young's inequality,
\[
|\delta_t|^{p-2}|u-a|^2\les M^{\beta(p-2)}\Xi^{p-2}(|u|^2+\Xi^2)\les M^{\beta(p-2)}(|u|^p+\Xi^p),
\]
and
\[
(|u-b|+|a-b|+|\delta_t|)^{p-2}|\delta_t|^2\les M^{2\beta}(|u|+|b|+\Xi)^{p-2}\Xi^2\les M^{2\beta}\big(|u|^p+\Xi^p+|b|^{p-2}\Xi^2\big).
\]
With $\eta(M)=M^{2\beta}+M^{\beta(p-2)}$ we get 
\[W(u-b,a_0-b)\les W(u-b,a-b)+\eta(M)(V_b^2(u)+V_b^2(\Xi)).\] 
Integrating in $t\in[0,1]$ and in $\pi$ over $\#_1$, and using $\pi(\#_1)\les1$,
\begin{equation}\label{lhs}
\int_{\#_1}W(y-x-b,\xi_r(x)-b)\,d\pi\les\int_\Omega\int_{t_-}^{t_+} W(\dot X-b,\xi_r(X)-b)\,dt\,d\pi+\eta(M)\hE .
\end{equation}

\step 6 {Conclusion.} Inserting \eqref{eq:T1}--\eqref{eq:T4} into \eqref{eq:estimateMain} and using \eqref{lhs},
\begin{align*}
\int_{\#_1}W(y-x-b,\xi_r(x)-b)\,d\pi\le{}& C\big(\delta+r^{\frac{p'-1}2}+D^{1/p}\big)\hE+C_{\delta,r}\big(M+M^\beta+M^2+D^{\beta/p}\big)\hE\\
&+C\eta(M)\hE+C_{\delta,r}D .
\end{align*}
We first choose $r$ so that $r^{(p'-1)/2}\le\delta$. Since $M\to0$ as $\eps\to0$ by \eqref{eq:LinftyBound}, we then choose $\eps$ so 
small that $D^{1/p}\le\delta$ and $C_{\delta,r}(M+M^\beta+M^2+D^{\beta/p})+C\eta(M)\le\delta$. This yields \eqref{conclharmapprox} after redefining $\delta$.
\end{proof}

\subsection{One-step improvement}\label{sec:onestep}

\begin{proposition}\label{onestep}
Let $C_0$ be the constant of Theorem \ref{harmapprox}. There exist $\theta\in(0,\frac18]$, $\eps_1>0$ and $C_1<\infty$ with the following property. 
Let $\pi$ be optimal for $c_b$, $R>0$ and $R'\in[4R,5R]$, and assume
 \begin{equation}\label{hypsmallonestep}
 \frac{|b|}R+\frac{1}{R^p} \lt(E_{b,4R}(\pi)+ D_{b,R'}\rt)\le \eps_1 .
 \end{equation}
Then there exists $\tilde{b}\in \R^d$ with
\begin{equation}\label{estimshift}
 V_b^2(\tilde{b})\le C_0\big(E_{b,4R}(\pi)+ D_{b,R'}\big),\qquad |\tilde b|\le\theta R,
\end{equation}
such that the image $\tilde\pi$ of $\pi$ under $(x,y)\mapsto(\tilde x,\tilde y)=(x-\tilde b,y)$,  satisfies
\begin{equation}\label{conclusiononestep}
 E_{b+ \tilde{b},\theta R}(\tilde{\pi})\le \frac12\, E_{b,4R}(\pi) + C_1 D_{b,R'}.
\end{equation}
\end{proposition}

\begin{proof}
By a scaling argument, we may assume $R=1$. Let $\delta\in(0,1)$ to be fixed below in terms of $\theta$, and let $\eps_1\le\eps(\delta)$, 
where $\eps(\delta)$ is given by Theorem \ref{harmapprox}. Let $\xi$ be the field given by Theorem \ref{harmapprox}, $\Xi=\|\xi\|_{C^{0,\beta}(B_{7/4})}$, 
so that $V_b^2(\Xi)\le C_0\hE$, and set $\tilde b=-\xi(0)$. 
Then $V_b^2(\tilde b)\le V_b^2(\Xi)\le C_0\hE$, which is the first part of \eqref{estimshift}. Since $|\tilde b|^p\le V_b^2(\tilde b)\le C_0\eps_1$, we have $|\tilde b|\le\theta$ if $\eps_1\le\theta^p/C_0$.

Let $b'=b+\tilde b=b-\xi(0)$ and, for a pair $(x,y)$, $u=y-x$ and $\tilde u=\tilde y-\tilde x=u-\xi(0)$. Since $|\tilde u|\le|u-b|+|b-\xi(0)|$ and $|b'|=|\xi(0)-b|$,
\begin{equation}\label{eq:newexcess}
V_{b'}^2(\tilde u)\le (|b'|+|\tilde u|)^{p-2}|\tilde u|^2\les\,W(u-b,\xi(0)-b).
\end{equation}
 $(\tilde x,\tilde y)\in\#_\theta\cap \spt\tilde \pi$ if $\tilde x\in B_\theta$, i.e.\ $x\in B_\theta(\tilde b)\subset B_{2\theta}$, or if $y\in B_\theta$. Hence, if $(\tilde x, \tilde y)\in \#_\theta$, then $(x,y)\in \#_{2\theta}$. By Lemma \ref{lem:Linfty}, $(x,y)\in\#_{2\theta}\cap\spt\pi$ satisfy $|x|\le2\theta+M\le3\theta$ provided $\eps_1$ is so small that $M\le\theta$; this uses that $\#_\theta$ also contains the pairs whose \emph{target} lies in the small ball.

Fix $(x,y)\in \#_{2\theta}\cap \spt\pi$ such a pair, let $a=\xi(x)$ and $\delta_0=\xi(0)-\xi(x)$, so $|\delta_0|\le\Xi(3\theta)^\beta$ and $|a|\le\Xi$. As in Step 5 of the proof of Theorem \ref{harmapprox},
\[
W(u-b,\xi(0)-b)\les W(u-b,a-b)+|\delta_0|^{p-2}|u-a|^2+(|u-b|+|a-b|+|\delta_0|)^{p-2}|\delta_0|^2 .
\]
We use $|u-a|^p\le W(u-b,a-b)$ (because $|u-b|+|a-b|\ge|u-a|$), $|u-b|\le|u-a|+|a-b|$ and Young's inequality $X^{p-2}Y^2\le X^p+Y^p$:
\begin{align*}
&|\delta_0|^{p-2}|u-a|^2+(|u-b|+|a-b|+|\delta_0|)^{p-2}|\delta_0|^2\\
&\qquad\les W(u-b,a-b)+|\delta_0|^p+(|b|+2\Xi)^{p-2}|\delta_0|^2\les W(u-b,a-b)+\theta^{2\beta}V_b^2(\Xi),
\end{align*}
where in the last step we used $\theta\le1$ and $\beta p\ge2\beta$. Consequently, for $(x,y)\in\#_{2\theta}\cap\spt\pi$,
\begin{equation}\label{eq:holdererror}
W(u-b,\xi(0)-b)\le C\,W(u-b,\xi(x)-b)+C\theta^{2\beta}\,V_b^2(\Xi).
\end{equation}
For $(x,y)\in\#_{2\theta}\cap\spt\pi$, we have $x\in B_{3\theta}$, so $\pi(\#_{2\theta})\le\lambda(B_{3\theta})\le C_m(3\theta)^d$ by Lemma \ref{lem:mass}(2), provided $(3\theta)^{p(d+1)}\ge5^d\Dat$, which holds for $\eps_1$ small depending on $\theta$. Integrating \eqref{eq:newexcess}--\eqref{eq:holdererror} and using \eqref{conclharmapprox} and \eqref{eq:estimxithm},
\[
E_{b',\theta}(\tilde\pi)\le C\theta^{-d}\int_{\#_{2\theta}}W(u-b,\xi(0)-b)\,d\pi\le C\theta^{-d}\big(\delta E+C_\delta D\big)+C_\ast\theta^{2\beta}C_0\,\hE .
\]
We now choose $\theta\le\frac18$ with $C_\ast C_0\theta^{2\beta}\le\frac14$, then $\delta=\delta(\theta)$ with $C\theta^{-d}\delta\le\frac14$, and finally $\eps_1$ small depending on $\theta$ and $\delta$ as required above. This gives \eqref{conclusiononestep} with $C_1=C\theta^{-d}C_\delta+\frac14$.
\end{proof}

\subsection{Iteration}\label{sec:iteration}

We now iterate Proposition \ref{onestep} to prove Theorem \ref{thm:iteration}, stated in the introduction. Recall that $E_{0,4R_0}(\pi)$ is the left-hand side of (H2).

\begin{proof}[Proof of Theorem \ref{thm:iteration}]
Constants $C$ depend only on $C_D,C_E$; $M_0$ is chosen at the very end. Let $\theta$ be as in Proposition \ref{onestep}, $\sigma=\theta/4$, $R_k=\sigma^kR_0$,
 and let $K\ge0$ be the largest integer such that
\begin{equation}\label{eq:stoprule}
R_k\ge M_0\,r_\ast\,(1+k)\qquad\text{for all }0\le k\le K .
\end{equation}
Since $R_0\ge M_0r_\ast$ and $R_k\to0$, $K$ is well defined; it depends only on $r_\ast,R_0,M_0$ and $\sigma$. 
We construct vectors $\bar b_0=0,\bar b_1,\dots,\bar b_K$ and denote by $\pi_k$ the image of $\pi$ under $(x,y)\mapsto(x-\bar b_k,y)$.

\step 1 {Data on shifted balls.} Let $0\le k\le K$ and assume $|\bar b_k|\le R_k$. Write $\lambda_k=\lambda(\cdot+\bar b_k)$, so that the marginals of 
$\pi_k$ are $\lambda_k$ and $\mu$. Let $\hat R_k$ be the dyadic number in $[9R_k,18R_k)$; then $\hat R_k\in[r_\ast,20R_0]$ since $R_k\ge M_0r_\ast$ by \eqref{eq:stoprule}.
 With $\rho=4R_k$ and $z\in\{0,\bar b_k\}$, we have $B_{2\rho}(z)\subset B_{9R_k}\subset Q_{\hat R_k}$ and $\rho\ge \hat R_k/16$.
  By (H1), $\hat R_k^{-p}\Dat^\square_{\hat R_k}\le C_DM_0^{-p}\le\eps_{\rm loc}$ for $M_0$ large, so Lemma \ref{inp:loc}, applied to $\lambda$ with centre $\bar b_k$ and to $\mu$ with centre $0$, yields
\[
\frac1{\rho}\int_{\rho}^{2\rho}\big(\Dat_{r}(\lambda_k)+\Dat_{r}(\mu)\big)dr\le 2C_{\rm loc}C_Dr_\ast^p .
\]
Hence there is a common good radius $R'_k\in[4R_k,5R_k]$ with $\Dat_{R'_k}(\lambda_k)+\Dat_{R'_k}(\mu)\le 16C_{\rm loc}C_Dr_\ast^p$. Thus, with $D_k:=D_{\bar b_k,R'_k}$ computed for $(\lambda_k,\mu)$,
\begin{equation}\label{eq:Dk}
D_k\le C_2\big(r_\ast^p+r_\ast^2|\bar b_k|^{p-2}\big).
\end{equation}

\step 2 {Induction.} Let $E_k=E_{\bar b_k,4R_k}(\pi_k)$ and $H_k=\max\{r_\ast,|\bar b_0|,\dots,|\bar b_k|\}$. Since $r_\ast\le H_k$, \eqref{eq:Dk} gives $D_k\le2C_2r_\ast^2H_k^{p-2}$. Set
\[
A=\max(C_E,4C_1C_2),\qquad C_3=\big(2^{\frac{p-2}2}C_0(A+2C_2)\big)^{1/2},\qquad C_4=1+C_3 .
\]
We prove by induction on $k\in\{0,\dots,K\}$ that
\begin{equation}\label{eq:induction}
E_k\le A\,r_\ast^2H_k^{p-2},\qquad H_k\le r_\ast(1+C_3k).
\end{equation}
For $k=0$, $\bar b_0=0$, $H_0=r_\ast$, and (H2) gives $E_0\le C_Er_\ast^p\le Ar_\ast^2H_0^{p-2}$. Assume \eqref{eq:induction} holds for some $k<K$. Then $|\bar b_k|\le H_k\le C_4r_\ast(1+k)\le C_4R_k/M_0\le R_k$, so Step 1 applies, and
\[
\frac{|\bar b_k|}{R_k}+\frac{E_k+D_k}{R_k^p}\le \frac{C_4}{M_0}+\frac{(A+2C_2)r_\ast^pC_4^{p-2}(1+k)^{p-2}}{M_0^pa^p(1+k)^p}\le\frac{C_4}{M_0}+\frac{(A+2C_2)C_4^{p-2}}{M_0^p}\le\eps_1
\]
for $M_0$ large, where we used \eqref{eq:stoprule}. Proposition \ref{onestep}, applied to $\pi_k$ with $b=\bar b_k$, $R=R_k$ and $R'=R'_k\in[4R_k,5R_k]$, gives $b_{k+1}$ with
\begin{equation}\label{eq:oneStepb}
V_{\bar b_k}^2(b_{k+1})\le C_0(E_k+D_k)\le C_0(A+2C_2)\,r_\ast^2H_k^{p-2}
\end{equation}
and, since $\theta R_k=4R_{k+1}$ and the image of $\pi_k$ under the shift by $b_{k+1}$ is $\pi_{k+1}$ with $\bar b_{k+1}=\bar b_k+b_{k+1}$,
\begin{equation}\label{eq:oneStepEnergy}
E_{k+1}\le\tfrac12E_k+C_1D_k\le\tfrac12Ar_\ast^2H_k^{p-2}+2C_1C_2r_\ast^2H_k^{p-2}\le Ar_\ast^2H_k^{p-2}\le Ar_\ast^2H_{k+1}^{p-2}.
\end{equation}
For the shift, if $H_{k+1}=H_k$ there is nothing to prove. Otherwise $H_{k+1}=|\bar b_k+b_{k+1}|\le|\bar b_k|+|b_{k+1}|$, and since $V_{\bar b_k}^2(b_{k+1})\ge2^{-\frac{p-2}2}(|\bar b_k|+|b_{k+1}|)^{p-2}|b_{k+1}|^2$, \eqref{eq:oneStepb} gives
\[
H_{k+1}^{p-2}|b_{k+1}|^2\le 2^{\frac{p-2}2}C_0(A+2C_2)\,r_\ast^2H_k^{p-2}\le C_3^2r_\ast^2H_{k+1}^{p-2}.
\]
Hence $|b_{k+1}|\le C_3r_\ast$ and $H_{k+1}\le H_k+C_3r_\ast\le r_\ast(1+C_3(k+1))$. This proves \eqref{eq:induction}; in particular $|\bar b_k|\le C_4r_\ast(1+k)$ for $0\le k\le K$.

\step 3 {Stopping radius and conclusion.} Since $R_K\ge M_0r_\ast$, we have $\sigma^KR_0\ge M_0r_\ast$, hence $1+K\le1+\log(R_0/r_\ast)/\log(1/\sigma)\les\log(R_0/r_\ast)$. By the maximality of $K$, $\sigma^{K+1}R_0=R_{K+1}<M_0r_\ast(K+2)$, so that
\[
R_K\le\sigma^{-1}M_0r_\ast(K+2)\les r_\ast\log(R_0/r_\ast).
\]
Conversely, taking logarithms in $\sigma^{K+1}R_0<M_0r_\ast(K+2)$ and using $\log(K+2)\le K+1$,
\[
\log(R_0/r_\ast)\le\log(3M_0)+(K+1)\big(1+\log(1/\sigma)\big)\le(1+K)\big(1+\log(1/\sigma)+\log(3M_0)\big),
\]
so that, by \eqref{eq:stoprule}, $R_K\ge M_0r_\ast(1+K)\ge c\,r_\ast\log(R_0/r_\ast)$ with $c>0$ depending only on $\sigma$ and $M_0$. At step $K$ the dimensionless smallness verified in Step 2 holds (with $\eps_1\le\eps_\infty$), and Lemma \ref{lem:Linfty},
 applied after rescaling by $R_K$ to $\pi_K$ with $b=\bar b_K$, gives $|\tilde y-\tilde x|\le MR_K\le R_K$ for all $(\tilde x,\tilde y)\in\spt\pi_K\cap\#_{7R_K/2}$. Let $(x,y)\in\spt\pi\cap\#_{R_K/2}$. If $x\in B_{R_K/2}$, then $\tilde x=x-\bar b_K\in B_{R_K}$ because $|\bar b_K|\le C_4r_\ast(1+K)\le R_K/2$ for $M_0\ge2C_4$; if $y\in B_{R_K/2}$, then $\tilde y=y\in B_{R_K/2}$. Consequently $(\tilde x,\tilde y)\in\#_{7R_K/2}$ in the shifted coordinates, and the residual displacement is $\tilde y-\tilde x=y-x+\bar b_K$. Therefore
\[
|y-x|\le|y-x+\bar b_K|+|\bar b_K|\le R_K+C_4r_\ast(1+K)\le C\,r_\ast\log(R_0/r_\ast).
\]
Choosing $C_\star\ge\max(C,2/c)$, we have $R_\star\le R_K/2$, hence $\#_{R_\star}\subset\#_{R_K/2}$, and \eqref{eq:finaldisplacement} follows. From now on $C_\star$ is fixed; by the bound on $R_K$ above, $R_\star\ge c_\star R_K$ with $c_\star>0$ depending only on $d,p,C_D,C_E$.

\step 4 {Excess at intermediate scales.} Let $0\le k\le K$ and $b_k=-\bar b_k$. Since $|\bar b_k|\le C_4r_\ast(1+k)\le R_k$, every pair $(x,y)\in\#_{3R_k}$ satisfies $x-\bar b_k\in B_{4R_k}$ or $y\in B_{4R_k}$. Thus, in the coordinates of $\pi_k$, it belongs to $\#_{4R_k}$ and its displacement is $y-x+\bar b_k=y-x-b_k$. Using $V_{\bar b_k}^2(z)\ge|z|^p$, \eqref{eq:induction} and $H_k\le C_4r_\ast(1+k)$, we obtain
\begin{equation}\label{eq:excessk}
\int_{\#_{3R_k}}|y-x-b_k|^p\,d\pi\le(4R_k)^dE_k\le C R_k^d\,r_\ast^p(1+k)^{p-2}.
\end{equation}
Let now $R\in[R_\star,R_0]$ and let $k$ be the largest integer in $\{0,\dots,K\}$ such that $R\le3R_k$. If $k<K$, then $R>3R_{k+1}\ge\sigma R_k$. If $k=K$, then $R\ge R_\star\ge c_\star R_K$. In both cases $R\ge cR_k$, and since $\#_R\subset\#_{3R_k}$, \eqref{eq:excessk} gives
\[
\Big(\frac1{R^d}\int_{\#_R}|y-x-b_k|^p\,d\pi\Big)^{\frac1p}\le C r_\ast(1+k)^{1-\frac2p}.
\]
Finally, $R\le3R_k=3\sigma^kR_0$ gives $k\log(1/\sigma)\le\log3+\log(R_0/R)$, hence $1+k\le C(1+\log(R_0/R))$, and $|b_k|\le C_4r_\ast(1+k)\le Cr_\ast(1+\log(R_0/R))$. This proves \eqref{eq:excessdecay} for a suitable $C_{\rm exc}$.
\end{proof}

\section{Higher integrability of the matching cost}\label{sec:higher}

Let $L\ge2$ with $n=L^d\in\N$, and equip $\T_L=(\R/L\Z)^d$ with the distance $d_{\T_L}$. Throughout the stochastic part of the paper we consider 
the binomial model: $\lambda=\sum_{i=1}^n\delta_{X_i}$ and $\mu=\sum_{j=1}^n\delta_{Y_j}$, where $(X_i)$ and $(Y_j)$ are independent families of
 i.i.d.\ uniform points on $\T_L$. Let $\pi_p$ be an optimal coupling for $W_p(\lambda,\mu)$ with cost $d_{\T_L}(x,y)^p$.

\begin{proof}[Proof of Theorem \ref{thm:main}]
In this proof $\les$ indicates a constant that depends only on $d,p,q$. Since  $d_{\T_L}\le\sqrt dL/2$ we may assume that $L\gg 1$. 

\step 1 {Stationarity.} By stationarity of the laws of $\lambda$ and $\mu$, 
for every deterministic $R<L/4$, averaging over translations and using Fubini gives
\begin{equation}\label{eq:stationaryreduction}
\EE\Big[\frac1{L^d}\int_{\T_L\times\T_L}d_{\T_L}(x,y)^q\,d\pi_p\Big]
=\frac1{|B_R|}\EE\Big[\int_{B_R\times\T_L}d_{\T_L}(x,y)^q\,d\pi_p\Big]
\le\EE[\kappa_{\lambda,R}D_R^q],
\end{equation}
where $\kappa_{\lambda,R}=\lambda(B_R)/|B_R|$ is as in \eqref{eq:kappa} and
\[
D_R=\sup\{d_{\T_L}(x,y):(x,y)\in\spt\pi_p,\ x\in B_R\},
\]
with the supremum of the empty set taken to be zero.  Since  $\lambda(B_R)$ has law $\operatorname{Bin}(n,|B_R|/n)$, 
\begin{equation}\label{eq:localmassmoments}
\EE \kappa_{\lambda,R}=1,\qquad
\EE \kappa_{\lambda,R}^2\les 1.
\end{equation}

\step 2 {The good event at the origin.} Set $A=2q+1$, let $r_L=(C_A\log L)^{1/d}$ be given by Theorem \ref{thm:momentlog}, and choose
\[
\ell=\frac L{100},\qquad R_0=\frac{\ell}{40\sqrt d},\qquad
R_L=C_\star^{-1}r_L\log(R_0/r_L),
\]
where $C_\star=C_\star(d,p,q)$ is the constant of Theorem \ref{thm:iteration} for the constants $C_D,C_E$ obtained in Step 4 below.
 For $L\gg1$, $|B_{R_L}|\gg1$ and $R_L\ll L$. Let $\mathcal G$ be the intersection of the events of Theorem \ref{thm:momentlog} 
 for $\lambda$ and $\mu$. Then 
 \[\PP(\mathcal G^c)\le2L^{-A}.\] 

On $\mathcal G$, Theorem \ref{thm:momentlog}(a) gives for every dyadic $R\in[r_L,L/2)$
\begin{equation}\label{eq:torusdata}
\Dat^\square_R(\lambda)+\Dat^\square_R(\mu)\les r_L^p.
\end{equation}
 By the triangle inequality, Theorem \ref{thm:momentlog}(c) also gives
\begin{equation}\label{eq:torusenergy}
L^{-d}\int d_{\T_L}^p\,d\pi_p
\les L^{-d}W_{p,\T_L}^p(\lambda,dx)
+L^{-d}W_{p,\T_L}^p(\mu,dx)\les r_L^p.
\end{equation}
By a straightforward adaptation of Lemma \ref{lem:Linfty}(2) to the torus, we have 
\begin{equation}\label{eq:aprioriLinfty}
|x-y|=d_{\T_L}(x,y)\le\ell\qquad
\text{for }(x,y)\in\spt\pi_p\text{ with }x\in B_{2\ell}\text{ or }y\in B_{2\ell}.
\end{equation}
In particular, considering the restriction of $\pi_p$ to $(B_{2\ell}\times \T_L)\cup(\T_L\times B_{2\ell})$, we can identify it with an Euclidean transport map between its marginals.

\step 3  {Iteration.} Since $Q_{20R_0}\subset B_{\ell}$, \eqref{eq:torusdata} gives (H1) of Theorem \ref{thm:iteration} with $r_\ast=r_L$ for $\pi_{p}$, 
at every dyadic scale $R\in[r_L,20R_0]$, with $C_D=C_D(d,p,q)$. By \eqref{eq:torusenergy} we also have 
\[
E_{0,4 R_0}(\pi_p) \les L^{-d} \int d_{\T_L}^p\,d\pi_p \les r_L^p,
\]
so (H2) holds with $C_E=C_E(d,p,q)$.
Moreover,
 for $L\gg1$ we also have $R_0\ge M_0r_L$. Since $R_L=R_\star$  the displacement estimate \eqref{eq:finaldisplacement} of 
  Theorem \ref{thm:iteration} gives, on $\mathcal G$,
\begin{equation}\label{eq:localdisplacement}
D_{R_L}\le C_\star\,r_L\log(R_0/r_L)
=C_\star^2R_L\le C(\log L)^{1+1/d}.
\end{equation}

\step 4 {Conclusion.} Apply \eqref{eq:stationaryreduction} with $R=R_L$. On $\mathcal G$ use \eqref{eq:localdisplacement}, 
and on its complement use $D_{R_L}\les L$. By \eqref{eq:localmassmoments} and Cauchy--Schwarz,
\begin{align*}
\EE\Big[\frac1{L^d}\int d_{\T_L}^q\,d\pi_p\Big]
&\les R_L^q\,\EE \kappa_{\lambda,R_L}+L^q\,\EE[\kappa_{\lambda,R_L}\one_{\mathcal G^c}]\\
&\les R_L^q+L^q(\EE \kappa_{\lambda,R_L}^2)^{1/2}\PP(\mathcal G^c)^{1/2}\\
&\les R_L^q+L^{q-A/2}.
\end{align*}
Since $A>2q$, this proves \eqref{mainclaim}.
\end{proof}

\section{The concentration result}\label{sec:concentration}

We now deduce a concentration estimate from the moment bound of Theorem \ref{thm:main}. We mostly follow the proof from \cite{d2026concentration} but include some details for the reader's convenience.
  Let $d\ge3$, $p\ge2$, $L\ge2$ with $n=L^d\in\N$. Recall that $F_{L,p}=L^{-d}\,W_{p,\T_L}^p(\lambda,\mu)$ and set $\gamma_{d,p}=2(p-1)\big(1+\frac1d\big)$, the exponent appearing in Theorem \ref{thm:concentration}.

\begin{proof}[Proof of Theorem \ref{thm:concentration}]
  We regard $F_{L,p}$ as a function of $(X_1,\dots,X_n,Y_1,\dots,Y_n)\in(\T_L)^{2n}$ and let $\sigma_p$ be the optimal
 permutation corresponding to the optimal matching $\pi_p$ (both being unique a.e.). It is not hard to prove that $F_{L,p}$ is a Lipschitz function of its variable and that at a.e. point cloud, 
  
\begin{align}\label{eq:pointcloud}
|\nabla F_{L,p}|^2=2p^2L^{-2d}\sum_{i=1}^nd_{\T_L}(X_i,Y_{\sigma_p(i)})^{2p-2}.
\end{align}

The uniform probability measure on $\T_L$ satisfies the Poincar\'e inequality with constant $L^2/(4\pi^2)$. This constant is preserved by tensorization, so it holds for the product measure on $(\T_L)^{2n}$ and the Lipschitz function $F_{L,p}$. In combination with \eqref{eq:pointcloud}, this implies,
\[
\Var(F_{L,p})\le\frac{L^2}{4\pi^2}\EE|\nabla F_{L,p}|^2
\le C_pL^{2-2d}\,\EE\Big[\sum_{i=1}^nd_{\T_L}(X_i,Y_{\sigma_p(i)})^{2p-2}\Big].
\]

Apply Theorem \ref{thm:main} with $q=2p-2$ and the optimal plan $\pi_p=\sum_i\delta_{(X_i,Y_{\sigma_p(i)})}$.  Since $q(1+\frac1d)=\gamma_{d,p}$, it gives
\[
\EE\Big[\sum_{i=1}^nd_{\T_L}(X_i,Y_{\sigma_p(i)})^{2p-2}\Big]
=L^d\,\EE\Big[\frac1{L^d}\int d_{\T_L}^q\,d\pi_p\Big]\les L^d(\log L)^{\gamma_{d,p}}.
\]
Thus $\Var(F_{L,p})\les L^{2-d}(\log L)^{\gamma_{d,p}}$, and Chebyshev's inequality gives
\[
\PP\big[|F_{L,p}-\EE F_{L,p}|\ge t\big]\le\frac{\Var(F_{L,p})}{t^2}.\qedhere
\]
\end{proof}


\appendix

\section{Non-linear PDE estimates}\label{app:PDE}

Throughout this appendix $R\in[2,3]$, $B=B_R$, $\nu$ is the outer normal and $\sigma_{d-1}$ the surface measure on 
$\partial B$. Constants are uniform in $R\in[2,3]$ and in $b\in\R^d$. For a density $h$ on $\partial B$ we write
 $H=\int_{\partial B}V_b^2(h)$; since $V_b^2(z)\ge|z|^p$, $H<\infty$ implies $h\in L^p(\partial B)$.

\subsection{Setting and statement}
\begin{definition}\label{def:Pr}
For $r\in(0,1]$ fix a smooth function $\rho_r\ge0$ on $SO(d)$, supported in $\{Q\colon|Q-\Id|\le r\}$, with $\rho_r(Q^{-1})=\rho_r(Q)$ and $\int\rho_r\,dQ=1$ (Haar measure). For $h\in L^1(\partial B)$ let $P_rh(x)=\int_{SO(d)}\rho_r(Q)h(Qx)\,dQ$.
\end{definition}

The image of the Haar measure under $Q\mapsto Qx$ is the normalized surface measure on $\partial B$. This gives the 
following properties: $\int_{\partial B}P_rh=\int_{\partial B}h$; $\int_{\partial B}(P_rh)\psi=\int_{\partial B}h\,P_r\psi$; $\int_{\partial B}\Psi(P_rh)\le\int_{\partial B}\Psi(h)$ for every convex $\Psi$; and, writing $P_rh(Q_0x)-P_rh(x)=\int(\rho_r(QQ_0^{-1})-\rho_r(Q))h(Qx)\,dQ$,
\begin{equation}\label{eq:Prsmooth}
\|P_rh\|_{C^{0,1}(\partial B)}\le C_r\int_{\partial B}|h| .
\end{equation}

A weak solution of \eqref{eq:PDEmain} with datum $h\in L^p(\partial B)$ and $\omega=|B|^{-1}\int_{\partial B} h$ is a $\phi\in W^{1,p'}(B)$ with $\xi=\nabla c_b^\ast(\nabla\phi)$ and
\begin{equation}\label{eq:weakform}
\int_B\langle\xi,\nabla\zeta\rangle=\int_{\partial B}h\,\zeta-\omega\int_B\zeta\qquad\text{for all }\zeta\in W^{1,p'}(B).
\end{equation}
Such a $\phi$ exists and is unique up to constants: it minimizes $\int_Bc_b^\ast(\nabla\psi)-\int_{\partial B}h\psi+\omega\int_B\psi$, which is invariant under $\psi\mapsto\psi+c$, strictly convex in $\nabla\psi$, and coercive on mean-zero functions by \eqref{eq:energyBoundVstar}, the trace inequality and the Poincar\'e inequality. We normalize $\int_B\phi=0$. The objects $\phi_r$, $\xi_r$ are defined with $h$ replaced by $P_rh$; note that the constant $\omega$ is unchanged.

\begin{theorem}\label{thm:regularity}
Let $b\in\R^d$, $r\in(0,\frac1{12}]$, and let $h$ be a density on $\partial B$ with $H=\int_{\partial B}V_b^2(h)<\infty$. Then:
\begin{enumerate}
\item[(i)] $\int_BV_{b,\ast}^2(\nabla\phi)+\int_B V_b^2(\xi)\les H$, and the same holds for $\phi_r,\xi_r$;
\item[(ii)] for $\psi\in\{\phi,\phi_r\}$, $\big|\int_{\partial B}(h-P_rh)\psi\big|\les r^{\frac{p'-1}2}H$.
\end{enumerate}
There are $\beta\in(0,1)$ and $C$, depending only on $d,p$, and $C_r$ depending also on $r$, such that
\begin{enumerate}
\item[(iii)] $V_b^2\big(\|\xi_r\|_{C^{0,\beta}(B_{7/4})}\big)\le CH$;
\item[(iv)] $V_b^2\big(\|\xi_r\|_{C^{0,\beta}(\overline B)}\big)\le C_rH$.
\end{enumerate}
\end{theorem}

Here $\|f\|_{C^{0,\beta}(U)}=\sup_U|f|+[f]_{C^{0,\beta}(U)}$. 
\subsection{Poincar\'e, trace inequalities and mollification}
Write $\Phi=V_{b,\ast}^2$, viewed as a function of $t=|y|\ge0$. It is convex and increasing, $\Phi(\ell t)\le\ell^2\Phi(t)$ for $\ell\ge1$, $\Phi(\ell t)\le\ell^{p'}\Phi(t)$ for $\ell\le1$, $\Phi(t)\le t^{p'}$, and $\Phi'(t)\le2(|b|^{p-1}+t)^{p'-2}t$.

\begin{lemma}\label{lem:traceEstimate}
For $\psi\in W^{1,p'}(B)$ with $\psi_B=\fint_B\psi$,
\begin{align*}
\int_B\Phi(|\psi-\psi_B|)+\int_{\partial B} \Phi(|\psi-\psi_B|)\lesssim \int_{B} \Phi(|\nabla \psi|).
\end{align*}
\end{lemma}
\begin{proof}
Both sides are continuous under $W^{1,p'}$-convergence, since $|\Phi(s)-\Phi(t)|\les(s^{p'-1}+t^{p'-1})|s-t|$ and traces converge in $L^{p'}(\partial B)$; hence we may assume $\psi\in C^1(\overline B)$. For $x,y\in B$, $\psi(x)-\psi(y)=\int_0^1\langle\nabla\psi(y+t(x-y)),x-y\rangle dt$ and $|x-y|\le6$. By Jensen's inequality,
\[
\Phi(|\psi(x)-\psi_B|)\le\fint_B\Phi(|\psi(x)-\psi(y)|)\,dy\les\fint_B\int_0^1\Phi(|\nabla\psi(y+t(x-y))|)\,dt\,dy .
\]
We integrate in $x\in B$. For $t\ge\frac12$ we change variables $x\mapsto z=y+t(x-y)$, with Jacobian $t^{-d}\le2^d$, and for $t<\frac12$ we change $y\mapsto z$, with Jacobian $(1-t)^{-d}\le2^d$; in both cases $z\in B$ by convexity. This proves the bound for $\int_B\Phi(|\psi-\psi_B|)$. For the trace, let $v=\psi-\psi_B$, $\theta\in {\mathbb S}^{d-1}$ and $\varrho\in[R/2,R]$. Then
\begin{align*}
\Phi(|v(R\theta)|)&\le\Phi(|v(\varrho\theta)|)+\int_{R/2}^R\Phi'(|v(t\theta)|)|\nabla\psi(t\theta)|\,dt\\
&\stackrel{\eqref{eq:YoungsV2}}{\les}\Phi(|v(\varrho\theta)|)+\int_{R/2}^R\big(\Phi(|v(t\theta)|)+\Phi(|\nabla\psi(t\theta)|)\big)dt.
\end{align*}
 Averaging over $\varrho\in[R/2,R]$ and integrating over $\theta$ concludes the proof.
\end{proof}

\begin{lemma}\label{lem:rotationPoincare}
For $\psi\in W^{1,p'}(B)$ and $r\in(0,\frac1{12}]$,
\[
\int_{\partial B}\Phi(|\psi-P_r\psi|)\les r^{p'-1}\int_B\Phi(|\nabla\psi|).
\]
\end{lemma}
\begin{proof}
By density we may assume $\psi\in C^1(\overline B)$. Let $x\in\partial B$ and $Q=\exp(A)$ with $|Q-\Id|\le r$, $|A|\les r$, and let $\tau\in(r,2r)$. Connect $x$ to $Qx$ by the path $\gamma_{x,Q,\tau}$ which goes radially from $x$ to $(1-\tau)x$, then along $t\mapsto(1-\tau)\exp(tA)x$, $t\in[0,1]$, and radially from $(1-\tau)Qx$ to $Qx$. Its length $\ell$ satisfies $\ell\le 12r\le1$. Parametrizing by arc length, $|\psi(x)-\psi(Qx)|\le\int_0^\ell|\nabla\psi(\gamma(s))|ds$, so by Jensen's inequality and $\Phi(\ell t)\le\ell^{p'}\Phi(t)$,
\[
\Phi(|\psi(x)-\psi(Qx)|)\le\ell^{p'-1}\int_0^\ell\Phi(|\nabla\psi(\gamma(s))|)\,ds\le(12r)^{p'-1}\int_{\gamma_{x,Q,\tau}}\Phi(|\nabla\psi|).
\]
By Jensen's inequality, we deduce that  $\Phi(|\psi(x)-P_r\psi(x)|)\le\int\rho_r(Q)\Phi(|\psi(x)-\psi(Qx)|)dQ$. We average the previous bound over $\tau\in(r,2r)$ and integrate over $x\in\partial B$ and $Q$. The two radial pieces contribute $\les\int_{B\setminus B_{R/2}}\Phi(|\nabla\psi|)$, since $Q$ preserves the surface measure. For the angular piece, $|\frac d{dt}\exp(tA)x|\le Cr$ and, by rotation invariance,
\begin{align*}
&\frac1r\int_r^{2r}\int_{\partial B}\int_0^1\Phi\big(|\nabla\psi((1-\tau)\exp(tA)x)|\big)\,Cr\,dt\,d\sigma_{d-1}(x)\,d\tau\\
&\qquad=C\int_r^{2r}\int_{\partial B}\Phi(|\nabla\psi((1-\tau)x)|)\,d\sigma_{d-1}\,d\tau\le C\int_B\Phi(|\nabla\psi|).
\end{align*}
Collecting the three pieces gives the claim.
\end{proof}

\subsection{Energy estimates}
\begin{lemma}\label{lem:energyEstimate}
Theorem \ref{thm:regularity}(i) holds.
\end{lemma}
\begin{proof}
Testing \eqref{eq:weakform} with $\zeta=\phi$ (recall $\int_B\phi=0$) gives $\int_B\langle\xi,\nabla\phi\rangle=\int_{\partial B}h\phi$. By the equality case of Fenchel's inequality, $\langle\xi,\nabla\phi\rangle=c_b^\ast(\nabla\phi)+c_b(\xi)$. By \eqref{eq:YoungsV} (second form) and Lemma \ref{lem:traceEstimate}, for $\eps\in(0,1]$,
\[
\int_{\partial B}h\phi\les\eps\int_{\partial B}\Phi(|\phi|)+\eps^{1-p}\int_{\partial B}V_b^2(h)\les \eps\int_B\Phi(|\nabla\phi|)+\eps^{1-p}H .
\]
By \eqref{eq:energyBoundV} and \eqref{eq:energyBoundVstar}, choosing $\eps$ small and absorbing, we obtain (i) for $\phi$. For $\phi_r$ we use in addition $\int_{\partial B}V_b^2(P_rh)\le H$.
\end{proof}

\begin{lemma}\label{lem:mollificationControl}
Theorem \ref{thm:regularity}(ii) holds. Consequently
\[
\Big|\int_{B}\langle\nabla\phi,\xi-\xi_r\rangle\Big|+\Big|\int_{B}\langle\nabla\phi_r,\xi-\xi_r\rangle\Big|\les r^{\frac{p'-1}2}H .
\]
\end{lemma}
\begin{proof}
Let $\psi\in\{\phi,\phi_r\}$. By self-adjointness, $\int_{\partial B}(h-P_rh)\psi=\int_{\partial B}h(\psi-P_r\psi)$. By \eqref{eq:YoungsV} with parameter $\eps\in(0,1]$, Lemma \ref{lem:rotationPoincare} and part (i),
\[
\Big|\int_{\partial B}h(\psi-P_r\psi)\Big|\le\eps H+C\eps^{-1}\int_{\partial B}\Phi(|\psi-P_r\psi|)\le\eps H+C\eps^{-1}r^{p'-1}H .
\]
We choose $\eps=r^{(p'-1)/2}$. For the consequence, subtract the weak formulations \eqref{eq:weakform} for $\xi$ and $\xi_r$ (same $\omega$) and test with $\zeta=\psi$.
\end{proof}

\subsection{Reduction to Neumann $p$-Laplace and known results}\label{sec:regreduction}
By scaling we may assume that $b=1$.  We then set
\[
u=\phi-\langle b,x\rangle,\qquad\text{so that}\qquad \xi=\nabla c_b^\ast(\nabla\phi)=b+|\nabla u|^{p'-2}\nabla u.
\]
 Then $u$ is a weak solution of the Neumann problem for the $p'$-Laplacian
\begin{equation}\label{eq:pprimeLaplace}
\nabla\cdot(|\nabla u|^{p'-2}\nabla u)=\omega\ \text{ in }B,\qquad |\nabla u|^{p'-2}\nabla u\cdot\nu=h-\langle b,\nu\rangle\ \text{ on }\partial B .
\end{equation}


We shall use the following regularity results for $q=p'\in(1,2]$.  

\begin{theorem}[Interior $C^{1,\beta}$ estimate]\label{inp:interior}
There are $\beta_1\in(0,1)$ and $C$, depending on $d,q$, such that every weak solution $v\in W^{1,q}(B_\varrho(x_0))$ of $\nabla\cdot(|\nabla v|^{q-2}\nabla v)=\omega$ with $\omega\in\R$ satisfies
\[
\sup_{B_{\varrho/2}(x_0)}|\nabla v|+\varrho^{\beta_1}[\nabla v]_{C^{0,\beta_1}(B_{\varrho/2}(x_0))}\le C\Big(\big(\fint_{B_\varrho(x_0)}|\nabla v|^q\big)^{1/q}+(\varrho|\omega|)^{\frac1{q-1}}\Big).
\]
\end{theorem}
Here $\les$ indicates a constant depending only on $d,q$. 

\begin{proof} The theorem essentially is follows from results in \cite{Lieberman1992survey} and a scaling argument. 

Let
\[
T=\Big(\fint_{B_\varrho(x_0)}|\nabla v|^q\Big)^{1/q}+(\varrho|\omega|)^{\frac1{q-1}} .
\]
If $T=0$, then $\nabla v=0$ a.e.\ and there is nothing to prove. Otherwise set
\[
w(y)=\frac{v(x_0+\varrho y)-(v)_{B_\varrho(x_0)}}{\varrho T},\qquad y\in B_1 .
\]
Then $\fint_{B_1}|\nabla w|^q\le1$ and $\nabla\cdot(|\nabla w|^{q-2}\nabla w)=\tilde\omega:=\varrho\omega/T^{q-1}$, with $|\tilde\omega|\le1$ by the choice of $T$. Since $w$ has mean zero, the Poincar\'e inequality gives $\|w\|_{L^q(B_1)}\les 1$.

Apply \cite[Theorem 2.1]{Lieberman1992survey} on the balls $B_{1/8}(y_0)$, $y_0\in B_{3/4}$,
to find $\sup_{B_{1/16}(y_0)}w\les 1+\|w\|_{L^q(B_1)}$. The function $-w$ solves the same equation with $\tilde\omega$ replaced by $-\tilde\omega$. Hence $\sup_{B_{3/4}}|w|= M\les 1$.

Consequently, applying \cite[Theorem 3.2]{Lieberman1992survey} in $\Omega=B_{3/4}$, $\Omega'=B_{5/8}$, in the form for bounded weak solutions described in the paragraph following it (p.~303)
gives $\beta_1=\beta_1(d,q)$, such that $[\nabla w]_{C^{0,\beta_1}(B_{5/8})}\les 1$.  Since $\fint_{B_1}|\nabla w|^q\le1$ a simple interpolation estimate gives also 
$\sup_{B_{1/2}}|\nabla w|\les 1$. This concludes the proof.

\end{proof}

\begin{theorem}[Conormal boundary $C^{1,\beta}$ estimate]\label{inp:boundary}
Let $\alpha=\frac12$. Assume $\int_{\partial B}g=\omega|B|$. There are $\beta_1\in(0,1)$ and $C$, depending on $d,q$, such that every weak solution $v\in W^{1,q}(B)$ of $\nabla\cdot(|\nabla v|^{q-2}\nabla v)=\omega$ in $B$, $|\nabla v|^{q-2}\nabla v\cdot\nu=g$ on $\partial B$, satisfies
\[
\|\nabla v\|_{C^{0,\beta_1}(\overline B)}\le C\Big(\|\nabla v\|_{L^q(B)}+\|g\|_{C^{0,\alpha}(\partial B)}^{\frac1{q-1}}+|\omega|^{\frac1{q-1}}\Big).
\]
\end{theorem}
\begin{proof}
The result is a consequence of \cite[Theorem 1]{Liebermann} and a standard boundary $L^\infty-W^{1,p}$-estimate, see e.g. \cite[Chapter 7]{giusti2003direct}, in combination with a scaling argument. Alternatively, one can directly apply \cite[Corollary 1.5]{Antonini2027}.


In either case, given a general solution let $T=\|\nabla v\|_{L^q(B)}+\|g\|_{C^{0,1/2}(\partial B)}^{1/(q-1)}+|\omega|^{1/(q-1)}$. If $T=0$, then $\nabla v=0$. Otherwise $w=(v-(v)_B)/T$ is a weak solution with datum $g/T^{q-1}$ and forcing $\omega/T^{q-1}$, and due to the results mentioned above $\|\nabla w\|_{C^{0,\beta_1}(\overline B)}\le 1$ for some $\beta_1$ depending only on $d,q$.
\end{proof}

\begin{theorem}[Linear Schauder estimates]\label{inp:linear}
Let $\Lambda\ge1$, $\beta_1\in(0,\frac12)$, and let $A$ be a symmetric matrix field with $\Lambda^{-1}\Id\le A\le\Lambda\Id$ and $[A]_{C^{0,\beta_1}}\le\Lambda$ on the relevant domain. The constants below are $C=C(d,\beta_1,\Lambda)$, uniform in $\varrho\in[\frac18,1]$ and in $R\in[2,3]$.
\begin{enumerate}
\item[(a)] If $w\in W^{1,2}(B_\varrho(x_0))$, $\varrho\in[\frac18,1]$, is a weak solution of $\nabla\cdot(A\nabla w)=\omega$ in $B_\varrho(x_0)$ with $\omega\in\R$, then $\|\nabla w\|_{C^{0,\beta_1}(B_{\varrho/2}(x_0))}\le C(\|\nabla w\|_{L^2(B_\varrho(x_0))}+|\omega|)$.
\item[(b)] If $w\in W^{1,2}(B)$ is a weak solution of the conormal problem $\nabla\cdot(A\nabla w)=\omega$ in $B$, $A\nabla w\cdot\nu=k$ on $\partial B$
 with $k\in C^{0,\beta_1}(\partial B)$, then $\|\nabla w\|_{C^{0,\beta_1}(\overline B)}\le C(\|\nabla w\|_{L^2(B)}+\|k\|_{C^{0,\beta_1}(\partial B)}+|\omega|)$.
\end{enumerate}
\end{theorem}

\begin{proof}
In both parts we may replace $w$ by $w-(w)_U$, $U\in\{B_\varrho(x_0),B\}$, which does not change the equation, 
and use the Poincar\'e inequality $\|w-(w)_U\|_{L^2(U)}\le C\|\nabla w\|_{L^2(U)}$.

(a) Write $\omega=\nabla\cdot F$ with $F(x)=\omega(x-x_0)/d$, so that $\nabla\cdot(A\nabla w-F)=0$ in $B_\varrho(x_0)$ and 
$\|F\|_{C^{0,\beta_1}(B_\varrho(x_0))}\le C|\omega|$. The claim is the interior Schauder estimate for divergence-form equations with $C^{0,\beta_1}$ 
coefficients and $C^{0,\beta_1}$ right-hand side $\nabla\cdot F$, see \cite{Giaquinta}.

(b) This is well-known, but we could not find a precise reference. One can obtain the result using a straightening of boundary procedure and the much more general  \cite[Theorem 1.4]{DongJeonVita2024} with weight exponent $0$ and $k=1$ with the choice $\varphi(t)=t^{\beta_1}$ for $\beta_1<\frac 1 2$, which implies $\sigma(t)\leq C t^{\beta_1}$ and consequently the $C^{1,\sigma}$-bound of \cite[Theorem 1.4]{DongJeonVita2024} implies the claim.


To straighten the boundary, cover a neighbourhood $\{R-\frac18<|x|\le R\}$ of $\partial B$ by the sets $\Psi_j^{-1}(B_{1/2}^+)$, $j=1,\dots,N(d)$, where $\Psi_j$ are smooth diffeomorphisms of open sets $U_j$ onto $B_1$ with $\Psi_j(U_j\cap B)=B_1^+$ and $\Psi_j(U_j\cap\partial B)=B_1'$, whose $C^2$ norms and those of their inverses are bounded uniformly in $R\in[2,3]$. For one chart $\Psi$, set $\tilde w=w\circ\Psi^{-1}$, $J=|\det D\Psi^{-1}|$ and let $J_\partial$ be the density of the surface measure of $\partial B$ with respect to $dy'$. Changing variables in the weak formulation with $\zeta$ supported in $U$ gives
\[
\nabla\cdot(\widetilde A\nabla\tilde w)=\tilde f\ \text{ in }B_1^+,\qquad(\widetilde A\nabla\tilde w)_d=-\tilde k\ \text{ on }B_1',
\]
in the weak sense, with $\widetilde A=J\,(D\Psi\,A\,D\Psi^T)\circ\Psi^{-1}$, $\tilde f=J\omega$ and $\tilde k=J_\partial\,k\circ\Psi^{-1}$ (the outer normal of $\{y_d>0\}$ is $-e_d$). The matrix $\widetilde A$ is symmetric, with ellipticity and $C^{0,\beta_1}$ bounds depending only on $d,\Lambda$. Define
\[
F(y)=e_d\Big[-\tilde k(y')+\int_0^{y_d}\tilde f(y',s)\,ds\Big].
\]
Then $\nabla\cdot F=\tilde f$ and $F_d=-\tilde k$ on $B_1'$, so
\[
\nabla\cdot(\widetilde A\nabla\tilde w-F)=0\ \text{ in }B_1^+,\qquad(\widetilde A\nabla\tilde w-F)_d=0\ \text{ on }B_1'.
\]
\cite[Theorem 1.4]{DongJeonVita2024} now applies with $\mathbf g=-F$ and in combination with (a) and a simple covering argument gives (b).


\end{proof}



\subsection{Proof of Theorem \ref{thm:regularity}(iii)--(iv)}

We  note the following elementary fact. For $1<q\le2$ the map $z\mapsto|z|^{q-2}z$ is $(q-1)$-H\"older continuous, so for any domain $U$ and $v\in C^{1,\beta_1}(U)$,
\begin{equation}\label{Holdergeneral}
\big[|\nabla v|^{q-2}\nabla v\big]_{C^{0,\beta_1(q-1)}(U)}\les[\nabla v]_{C^{0,\beta_1}(U)}^{q-1}.
\end{equation}

\emph{Strategy.} The case $b=0$ being already covered by  \cite{koch2024geometric}, we may assume that $b\neq0$. By scaling we may further assume that  $|b|=1$. 
\begin{itemize}
  \item \emph{Moderate energy} If $H\ges 1$, the classical estimates for the $p'-$Laplacian give $\|\xi_r\|\les H^{1/p}$, which is the correct ($p$-th power) behaviour of $V_1^2$ away from zero.
\item \emph{Small energy} If $H\ll 1$, combining the energy bound and the regularity estimate for the $p'-$Laplacian we get that $|\nabla \phi_r|\ll 1$ so that the equation is uniformly ellliptic and the linear Schauder estimates apply.
\end{itemize}
The interior estimate (iii) and the global estimate (iv) follow the same scheme. The only differences are whether we rely on Theorem \ref{inp:interior} or Theorem \ref{inp:boundary} and whether constants may depend on $r$. 

\medskip
\emph{Notation and elementary bounds.} Let $\beta_1\in(0,\frac12)$ be a common exponent for Theorem \ref{inp:interior} and Theorem \ref{inp:boundary}. 
Theorem \ref{inp:linear} is used with this $\beta_1$. Set $\beta=\beta_1(p'-1)\le\beta_1$.  The goal is
\begin{equation}\label{eq:goalreg}
V_1^2(X_U)\le C_UH,\qquad X_U:=\|\xi_r\|_{C^{0,\beta}(U)},
\end{equation}
where $U=B_{7/4}$ with $C_U$ independent of $r$ for (iii), and $U=\overline B$ with $C_U=C_r$ for (iv). Four elementary bounds will be used repeatedly.
\begin{enumerate}
\item[(E1)] Since $|z|^{p'}\les V_{1,\ast}^2(z)+1$, part (i) gives
\begin{equation}\label{eq:Lpbound}
\int_B|\nabla u_r|^{p'}\les 1+H.
\end{equation}
\item[(E2)] For  every $z$, $V_1^2(z)=(1+|z|^2)^{\frac{p-2}2}|z|^2\ge|z|^2$ and, since $p\ge2$, $V_1^2(z)\ge|z|^p$. 
By H\"older's inequality, $\int_{\partial B}|h|\les(\int_{\partial B}|h|^2)^{1/2}\les H^{1/2}$ and $\int_{\partial B}|h|\les(\int_{\partial B}|h|^p)^{1/p}\les H^{1/p}$. Since $|\omega|\le|B|^{-1}\int_{\partial B}|h|$, we get
\begin{equation}\label{eq:L1bound}
|\omega|+\int_{\partial B}|h|\les\min\{H^{1/2},H^{1/p}\}.
\end{equation}
\item[(E3)] Let $U\in\{B_{15/8},\overline B\}$, $K\ge 1$ and  $G\colon U\to\R^d$ with $[G]_{C^{0,\beta_1}(U)}\le K$. For every $\eps>0$ there exists $\delta=\delta(K,\eps,\beta_1)$ such that if $\int_U V_{1,*}^2(G)\le \delta$ then $ \sup_U |G|\le \eps$. This follows from a simple interpolation argument. 
\item[(E4)] Let  let $U'\subset B$ be open and assume
\[
\sup_{U'}|\nabla\phi_r|\le\tfrac14,\qquad [\nabla\phi_r]_{C^{0,\beta_1}(U')}\le K .
\]
Then $\xi_r=A\nabla\phi_r$ in $U'$, where $A$ is a symmetric matrix field with $\Lambda^{-1}\Id\le A\le\Lambda\Id$ and $[A]_{C^{0,\beta_1}(U')}\le\Lambda$ for some $\Lambda=\Lambda(d,p,K)$. In particular $\nabla\cdot(A\nabla\phi_r)=\omega$ in $U'$, and if $U'=B$, then $A\nabla\phi_r\cdot\nu=P_rh$ on $\partial B$.
 Moreover $\|\xi_r\|_{C^{0,\beta_1}(U')}\les\|\nabla\phi_r\|_{C^{0,\beta_1}(U')}$. This result is obtained by a simple Taylor expansion of the gradient of the  cost around $b$.
\end{enumerate}

\step 1 {Interior estimate, moderate energy.} Let $H_0 \ll 1$ and assume $H\ge H_0$.
\begin{enumerate}
\item[(a)] By \eqref{eq:Lpbound}, $\|\nabla u_r\|_{L^{p'}(B)}\le C_{H_0}H^{1/p'}$, and by (E2),
\[
|\omega|^{p-1}\le C_{H_0}H^{(p-1)/p}=C_{H_0}H^{1/p'}.
\] Theorem \ref{inp:interior} (with $\frac1{q-1}=p-1$) and a simple covering argument give
\[
\|\nabla u_r\|_{C^{0,\beta_1}(B_{7/4})}\le C_{H_0}H^{1/p'}.
\]
\item[(b)] Since $\xi_r-b=|\nabla u_r|^{p'-2}\nabla u_r$, we have $\sup|\xi_r-b|\le(\sup|\nabla u_r|)^{p'-1}$, and \eqref{Holdergeneral} gives the H\"older seminorm. Since $(p'-1)/p'=1/p$, we obtain
\[
X_{B_{7/4}}\le|b|+C_{H_0}H^{1/p}\le C_{H_0}H^{1/p}.
\]
\item[(c)] Since $V_1^2(X)\les X^2+X^p$ and $H\ge H_0$, $p\ge2$, we get $V_1^2(X_{B_{7/4}})\le C_{H_0}(H^{2/p}+H)\le C_{H_0}H$.
\end{enumerate}

\step 2 {Interior estimate, small energy.} Let $H\le H_0\ll1$, where $H_0\le1$.
\begin{enumerate}
\item[(a)] \emph{A priori bound.} For $x_0\in B_{15/8}$ we have $B_{1/8}(x_0)\subset B$. By \eqref{eq:Lpbound} and \eqref{eq:L1bound} the right-hand side in Theorem \ref{inp:interior} (applied to $v=u_r$, $\varrho=\frac18$) is bounded by a constant. By a covering argument,
 $[\nabla\phi_r]_{C^{0,\beta_1}(B_{15/8})}=[\nabla u_r]_{C^{0,\beta_1}(B_{15/8})}\le K=K(d,p)$.
\item[(b)] \emph{Smallness.} By part (i) of Theorem \ref{thm:regularity}, $\int_{B_{15/8}}V_{1,\ast}^2(\nabla\phi_r)\les H$. By (E3) with $G=\nabla\phi_r$,
\[
\sup_{B_{15/8}}|\nabla\phi_r|\le\tfrac14
\]
provided $H_0$ is small enough.
\item[(c)] \emph{Linear estimate.} We may thus  apply (E4) on $U'=B_{15/8}$. Since $|\nabla\phi_r|\le1$ there, $|\nabla\phi_r|^2\les V_{1,\ast}^2(\nabla\phi_r)$, so $\|\nabla\phi_r\|_{L^2(B_{15/8})}\les H^{1/2}$. Theorem \ref{inp:linear}(a) on the balls $B_{1/8}(x_0)$, $x_0\in B_{7/4}$, together with \eqref{eq:L1bound} and covering, 
gives
\[
\|\nabla\phi_r\|_{C^{0,\beta_1}(B_{7/4})}\les\|\nabla\phi_r\|_{L^2(B_{15/8})}+|\omega|\les H^{1/2}.
\]
\item[(d)] \emph{Conclusion.} For $\|\nabla\phi_r\|_{C^{0,\beta_1}(B_{7/4})}\ll 1$, we have $\|\nabla\phi_r\|_{C^{0,\beta_1}(B_{7/4})}\simeq X_{B_{7/4}}$ to that  $X_{B_{7/4}}\les H^{1/2}\le1$. Since $V_1^2(X)\les X^2$ for $X\les1$, we get $V_1^2(X_{B_{7/4}})\les H$.
\end{enumerate}

In Steps 1 and 2, $h$ entered only through $H$ (via part (i) and \eqref{eq:L1bound}, and $\int V_1^2(P_rh)\le H$). 
The constants therefore do not depend on $r$. This proves (iii).

\step 3 {Global estimate.} We repeat Steps 1 and 2 on $U=\overline B$. The differences are as follows.
\begin{itemize}
\item Theorem \ref{inp:boundary} replaces Theorem \ref{inp:interior}. Its datum is the Neumann datum of $u_r$ in \eqref{eq:pprimeLaplace}, $g=P_rh-\langle b,\nu\rangle$, with forcing $\omega$; note that this pair satisfies the Neumann compatibility condition. By \eqref{eq:Prsmooth}, $\|g\|_{C^{0,\alpha}(\partial B)}\le C_r\int|h|+C$. This is $\le C_r$ if $H\le1$, and $\le C_rH^{1/p}$ if $H\ge H_0$ by (E2).
\item Theorem \ref{inp:linear}(b) replaces Theorem \ref{inp:linear}(a). By (E4) with $U'=B$, its conormal datum is $k=P_rh$, and by \eqref{eq:Prsmooth} and (E2), $\|k\|_{C^{0,\beta_1}(\partial B)}\le C_r\int|h|\le C_rH^{1/2}$ for $H\le1$.
\item Constants, and the threshold $H_0=H_0(r)$, may depend on $r$.
\end{itemize}

\section{Moment bound}\label{app:moments}

In this appendix $d\ge3$ and $p\ge1$. For a cube $Q$ and a measure $\nu$ we recall that 
 $\kappa_Q=\nu(Q)/|Q|$. For a cube $U$ of side $\ell_0$ and $i\ge0$, the cells of 
 level $i$ are the $2^{id}$ dyadic subcubes of $U$ of side
  $\ell_i=2^{-i}\ell_0$, and $\nu_i=\sum_{Q\text{ of level }i}\kappa_Q\one_Q\,dx$. We first recall the following  deterministic 
  bound which can be obtained for instance as a consequence of the Benamou-Brenier formula, see e.g. \cite[Lemma 3.4]{goldman2021convergence}.
\begin{lemma}\label{lem:refine}
Let $Q$ be a cube of side $\ell$ with children $Q_1,\dots,Q_{2^d}$, let $\kappa_1,\dots,\kappa_{2^d}\ge0$, $\beta=\sum_j\kappa_j\one_{Q_j}$ and $\kappa_Q=2^{-d}\sum_j\kappa_j$. Then
\begin{equation}\label{eq:refinegeneral}
W_p^p(\beta\,dx,\kappa_Q\one_Q\,dx)\les\ell^p|Q|\Big(\sum_j|\kappa_j-1|^p+|\kappa_Q-1|^p\Big),
\end{equation}
and if $\kappa_Q\ge\frac12$, then $W_p^p(\beta\,dx,\kappa_Q\one_Q\,dx)\les\ell^p|Q|\max_j|\kappa_j-\kappa_Q|^p$.
\end{lemma}

For the terminal step we use a trivial bound: if $\ell_I$ is the side of the cells of level $I$, then
\begin{equation}\label{eq:terminal}
W_p^p(\nu\rst U,\nu_I)\les\sum_{Q\text{ of level }I}\ell_I^p\nu(Q)=\ell_I^p\nu(U).
\end{equation}

For the probabilistic part we use Bernstein's inequality \cite{BouLuMa}: if $X$ is binomial with mean $v$, then for $t\ge0$
\begin{equation}\label{eq:bernstein}
\PP(|X-v|\ge t)\le2\exp\Big(-\frac{t^2}{2(v+t/3)}\Big).
\end{equation}
The same bound holds for the auxiliary variable $X\sim\operatorname{Pois}(v)$ used in the proof of Theorem \ref{thm:momentunit}.

\subsection{Elementary estimate with a logarithmic cutoff}
Let $\nu$ be $\lambda$ or $\mu$ of Section \ref{sec:higher}; thus $\nu(\T_L)=L^d$. For $R<L/2$, we identify the centred cube $Q_R$ with its image in $\T_L$. 

\begin{theorem}\label{thm:momentlog}
Let $A>0$. There are $C_A\ge1$ depending on $d,A$, $C$ depending on $d,p$, and $L_0$ depending on $d,p,A$, such that the following holds for $L\ge L_0$. Let $r_L=(C_A\log L)^{1/d}$, $\Lambda=(A+d+1)\log L$. With probability at least $1-L^{-A}$:
\begin{enumerate}
\item[(a)] $\Dat^\square_R(\nu)\le Cr_L^p$ for all dyadic $R\in[r_L,L/2)$;
\item[(b)] $\kappa_{Q_R}\in[\frac12,\frac32]$ for all $R$ as in (a);
\item[(c)] $L^{-d}W_{p,\T_L}^p(\nu,dx)\le Cr_L^p$.
\end{enumerate}
\end{theorem}
\begin{proof}
Let $\mathcal Q$ be the family of all cells of level $i\le I$ of the centred cubes $U=Q_R$, $R\in[r_L,L/2)$ dyadic, and of $U=[0,L)^d$, where for each $U$ the level $I$ is such that $\ell_I\in[r_L,2r_L)$. 
Each such partition has at most $C(R/r_L)^d$ cells over all its levels. Therefore $\#\mathcal Q\les(L/r_L)^d\le L^d$. For $Q\in\mathcal Q$, $\nu(Q)$ is binomial with parameters $L^d$ and $|Q|/L^d$, hence mean 
$|Q|$.  By \eqref{eq:bernstein} with $t=2\sqrt{|Q|\Lambda}+\Lambda$, which satisfies $t^2\ge2\Lambda(|Q|+t/3)$, and a union bound, with probability at least $1-2CL^de^{-\Lambda}\ge1-L^{-A}$ we have for every $Q\in\mathcal Q$ of side $\ell$
\[
|\kappa_Q-1|\le2\sqrt{\Lambda/\ell^d}+\Lambda/\ell^d\le\delta_\ell:=3\sqrt{\Lambda/\ell^d},
\]
where we choose $C_A=64(A+d+1)$, so that $\Lambda/\ell^d\le\frac1{64}$ and $\delta_\ell\le\frac38$ for $\ell\ge r_L$. This gives (b).

For (a), let $U=Q_R$, $\ell_0=2R$. By the triangle inequality for $W_p$ (the $\nu_i$ have the same mass),
\[
W_p(\nu\rst U,\kappa_U\one_Udx)\le W_p(\nu\rst U,\nu_I)+\sum_{i<I}W_p(\nu_{i+1},\nu_i).
\]
By subadditivity and Lemma \ref{lem:refine} (all densities lie in $[\frac12,\frac32]$), $W_p^p(\nu_{i+1},\nu_i)\les\ell_i^p|U|\delta_{\ell_{i+1}}^p$, hence
\[
\Big(\frac{W_p^p(\nu_{i+1},\nu_i)}{|U|}\Big)^{1/p}\les\ell_i\delta_{\ell_{i+1}}\les\sqrt\Lambda\,\ell_{i+1}^{1-\frac d2}.
\]
Since $d>2$, the sum over $i<I$ is geometric and dominated by its last term, so it is $\les\sqrt\Lambda r_L^{1-d/2}=r_L(\Lambda/r_L^d)^{1/2}\le r_L$. By \eqref{eq:terminal}, $(|U|^{-1}W_p^p(\nu\rst U,\nu_I))^{1/p}\les\ell_I\kappa_U^{1/p}\les r_L$. Minkowski's inequality gives $R^{-d}W_p^p(\nu\rst U,\kappa_U\one_Udx)\les r_L^p$. For the density term, $R^p|\kappa_U-1|^p\les R^p(\Lambda/R^d)^{p/2}\le\Lambda^{p/2}r_L^{p(1-d/2)}\le r_L^p$, again since $d>2$ and $R\ge r_L$. Part (c) is proved in the same way with $U=[0,L)^d$, since Euclidean transport in $[0,L)^d$ dominates transport on $\T_L$.
\end{proof}

\subsection{The unit-scale estimate}
 We now prove a stronger estimate using partitions down to unit scale for the binomial model. The main difficulty is that a simple union bound is not sufficient to treat sub-logarithmic scales which are treated by a CLT type argument.
We start with a simple tail estimate for sums of independent but not necessarily bounded random variables. This result is most likely already known.
\begin{lemma}\label{lem:truncation}
Let $Z_1,\dots,Z_N$ be independent with $\PP(|Z_j|\ge t)\le C_Ze^{-t/4}$ for all $t\ge0$, and let $m_p=\max_j\EE|Z_j|^p$. Then
\[
\PP\Big[\frac1N\sum_{j=1}^N|Z_j|^p>m_p+1\Big]\le Ce^{-cN^{1/(p+1)}},
\]
with $C,c$ depending on $p$ and $C_Z$.
\end{lemma}
\begin{proof}
Let $T=N^{p/(p+1)}$ and $W_j=\min(|Z_j|^p,T)$. Then
\[
\PP(\exists j\colon|Z_j|^p>T)\le C_ZNe^{-N^{1/(p+1)}/4}.
\] The variables $W_j-\EE W_j$ are independent, bounded by $T$ and have variance at most $\EE|Z_j|^{2p}\le C$. By Bernstein's inequality \cite{BouLuMa},
\[
\PP\Big[\sum_j(W_j-\EE W_j)\ge N\Big]\le\exp\Big(-\frac{N^2}{2(CN+TN/3)}\Big)\le\exp(-cN/T)=\exp(-cN^{1/(p+1)}).
\]
Outside both events, $\sum_j|Z_j|^p=\sum_jW_j\le N(m_p+1)$.
\end{proof}

\begin{theorem}\label{thm:momentunit}
Let $L\ge2$ with $n=L^d\in\N$, and let $\nu=\sum_{i=1}^n\delta_{X_i}$ for i.i.d.\ uniform points in  $\T_L$.
 Set 
\[
\eta=\frac{d(d-2)}{2d+(p+1)(d-2)}.
\]
There are $C,c>0$ depending only on $d,p$ such that for every cube $Q$ of side $R\in[1,L/2]$,
\[
\PP\Big[R^{-d}W_p^p(\nu\rst Q,\kappa_Q\one_Qdx)+R^p|\kappa_Q-1|^p>C\Big]\le Ce^{-cR^\eta}.
\]
\end{theorem}
\begin{proof}
\emph{Auxiliary Poisson estimate.} First let $\Pi$ be a unit-intensity Poisson process on $\T_L$. 
In this part of the proof, the densities $\kappa_{Q'}=\Pi(Q')/|Q'|$ and the piecewise constant measures $\Pi_i$ are defined using $\Pi$. Let $\ell_i=2^{-i}R$, $i=0,\dots,I$, where $I\ge0$ is such that $\ell_I\in[1,2)$,
 and 
 \[Y_\ell=N_\ell^{-1}\sum_{Q'}|\kappa_{Q'}-1|^p,\]
  the sum being over the $N_\ell=(R/\ell)^d$ cells of side $\ell$.

\emph{Deterministic bound.} Let $a=\frac{p(d-2)}8$. Put $x_i=(|Q|^{-1}W_p^p(\Pi_{i+1},\Pi_i))^{1/p}$ for $i<I$ and $x_I=(|Q|^{-1}W_p^p(\Pi,\Pi_I))^{1/p}$. By Minkowski's and, if $p>1$, H\"older's inequalities,
\[
\frac{W_p^p(\Pi,\kappa_Q\one_Qdx)}{|Q|}\le\Big(\sum_{i=0}^Ix_i\Big)^p\le\Big(\sum_{i}\ell_i^{-\frac a{p-1}}\Big)^{p-1}\sum_i\ell_i^ax_i^p\les\sum_i\ell_i^ax_i^p ;
\]
for $p=1$ the final bound $\sum_ix_i\le\sum_i\ell_i^ax_i$ is immediate since $\ell_i\ge1$.
By \eqref{eq:refinegeneral}, which requires no lower bound on the densities, and since each cell of side $\ell_{i+1}$ has
 exactly one parent, $x_i^p\les\ell_i^p(Y_{\ell_i}+Y_{\ell_{i+1}})$ for $i<I$; by \eqref{eq:terminal}, $x_I^p\les \kappa_Q$. Hence
\begin{equation}\label{eq:unitdet}
R^{-d}W_p^p(\Pi,\kappa_Q\one_Qdx)\les\kappa_Q+\sum_{i=0}^I\ell_i^{p+a}Y_{\ell_i}.
\end{equation}

\emph{Large scales.} These are treated as in the proof of Theorem \ref{thm:momentlog} by a union bound.  Let $\zeta=\frac{d-2}2$ and $\delta\in(0,1)$ to be fixed. For $\ell\ge R^\delta$ we apply \eqref{eq:bernstein} with $v=\ell^d$ and $t=\ell^{d-1-2a/p}\le v$: $\PP(|\kappa_{Q'}-1|>\ell^{-1-2a/p})\le2\exp(-\frac38\ell^{d-2-4a/p})=2\exp(-\frac38\ell^\zeta)$. On the event that this holds for all cells of all these levels, $Y_\ell\le\ell^{-p-2a}$; its complement has probability at most $2R^d(I+1)\exp(-\frac38R^{\delta\zeta})$.

\emph{Small scales.} For $\ell<R^\delta$ let $Z_{Q'}=(\Pi(Q')-\ell^d)\ell^{-d/2}$. At a fixed level these are independent (this is where it is important that $\Pi$ is a Poisson point process), and by \eqref{eq:bernstein}, for $t\ge1$, $\PP(|Z_{Q'}|\ge t)\le2\exp(-\min(\frac{t^2}4,\frac{3t\ell^{d/2}}4))\le2e^{-t/4}$, so the hypothesis of Lemma \ref{lem:truncation} holds with $C_Z=2e^{1/4}$, uniformly in $\ell\ge1$; also $m_p\le C$. Since $N_\ell\ge R^{d(1-\delta)}$, with probability at least $1-C(I+1)\exp(-cR^{d(1-\delta)/(p+1)})$ we have, for all these levels,
\[
Y_\ell=\ell^{-\frac{dp}2}\frac1{N_\ell}\sum_{Q'}|Z_{Q'}|^p\le(m_p+1)\ell^{-\frac{dp}2}\les \ell^{-p-2a},
\]
where we used $\frac{dp}2\ge p+2a$, i.e.\ $d\ge2$.

\emph{Top cell.} By \eqref{eq:bernstein} with $t=R^{d-1}$, $\PP(R|\kappa_Q-1|>1)\le2\exp(-cR^{d-2})$; on the complement $\kappa_Q\le2$.

\emph{Conclusion.} On the intersection of these events, \eqref{eq:unitdet} gives 
\[R^{-d}W_p^p\les \sum_i\ell_i^{-a}\les 1 \qquad \textrm{and } \qquad  R^p|\kappa_Q-1|^p\le1.\]
 Up to polynomial factors, which are absorbed by decreasing $c$, the failure probability is at most
\[
C\exp\Big(-cR^{\min(\delta\zeta,\,d(1-\delta)/(p+1),\,d-2)}\Big).
\] Choosing $\delta\zeta=\frac{d(1-\delta)}{p+1}$, i.e.\ $\delta=\frac{2d}{2d+(p+1)(d-2)}$, gives the exponent $\delta\zeta=\eta\le d-2$.

\emph{Transfer to the binomial process.} Let  $v=|Q|=R^d$. Since $|\T_L|=n$ and $R\le L/2$, we have $v\le2^{-d}n\le n/2$. Conditional on $\Pi(\T_L)=n$, $\Pi$ and   $\nu$ have the same law. For any event $A$ depending only on the restriction to $Q$, independence of $\Pi\rst Q$
 and $\Pi\rst(\T_L\setminus Q)$ gives
\begin{equation}\label{eq:localconditioning}
\PP_{\rm bin}(A)
=\EE_{\rm Pois}\left[\one_A\,
\frac{\PP(\operatorname{Pois}(n-v)=n-\Pi(Q))}
     {\PP(\operatorname{Pois}(n)=n)}\right],
\end{equation}
where probabilities of negative counts are zero. Stirling's formula gives $\PP(\operatorname{Pois}(n)=n)\ge cn^{-1/2}$, and the maximum point probability of $\operatorname{Pois}(n-v)$ is at most $C(n-v)^{-1/2}$. Hence, since $n-v\ge n/2$,
\[
\PP_{\rm bin}(A)\le C\sqrt{\frac{n}{n-v}}\,\PP_{\rm Pois}(A)
\le C\PP_{\rm Pois}(A).
\]
Apply this inequality to the final failure event of the auxiliary estimate, which depends only on $\Pi\rst Q$
 since all the cells involved are subcubes of $Q$. The binomial bound therefore has the same exponent $\eta$, with constants
  independent of $L$. 
\end{proof}

\begin{corollary}\label{cor:minimalradius}
Let $L\ge2$ with $n=L^d\in\N$, let $\nu=\sum_{i=1}^n\delta_{X_i}$ for i.i.d.\ uniform points in  $\T_L$, and let $C$ be the constant of Theorem \ref{thm:momentunit}. 
Define $r_{\ast,L}$ as the smallest dyadic $r\ge1$ such that for every dyadic $R\in[\max(2,r),L]$,
\[
R^{-d}W_p^p(\nu\rst Q_{R/2},\kappa_{Q_{R/2}}dx)+R^p|\kappa_{Q_{R/2}}-1|^p\le C.
\]
 There are $C',c'>0$, depending only on $d,p$, such that for $r\ge 1$,
\[
\PP(r_{\ast,L}>r)\le C'e^{-c'r^\eta}.
\]
Moreover, for every $\eps\in(0,c')$,
\[
\sup_{L\ge2:\,L^d\in\N}\EE\left[\exp(\eps r_{\ast,L}^\eta)\right]<\infty.
\]
\end{corollary}
\begin{proof}
If $r_{\ast,L}>r$, an estimate must fail at some dyadic $R\in[2,L]$ with $R>r/2$: otherwise the largest dyadic number at most $r$ would satisfy the defining condition. Theorem \ref{thm:momentunit} and a union bound therefore give
\[
\PP(r_{\ast,L}>r)
\le\sum_{\substack{2\le R\le L,\ R>r/2\\R\ \text{dyadic}}}Ce^{-cR^\eta}
\le C'e^{-c'r^\eta}.
\]
The constants are uniform in $L$. Integrating this tail bound for $r_{\ast,L}^\eta$ gives the stated uniform exponential moment for $\eps<c'$.
\end{proof}

\section*{Acknowledgments}
This research benefited from the support of the FMJH Program Gaspard Monge for optimization and operations research and their interactions with data science.

\section*{AI disclosure}
Most of the proofs in this paper were obtained without relying on AI tools. Large language models (LLMs) were instead used to streamline and simplify some of the proofs. In particular, in a first version of the paper we had reproved the regularity results for the geometrically linearized problem following the strategy of Duzaar and Mingione \cite{duzaar2004}, relying also on the approach for the Neumann problem in Irving and the second author \cite{KI2023}. LLMs were then used to help with the writing of the paper. The authors take full responsibility for the content of the paper.

\bibliographystyle{acm}
\IfFileExists{OT.bib}{\bibliography{OT}}{\bibliography{../OT}}
\end{document}